\documentclass{article}

\usepackage[cp1250]{inputenc}
\usepackage[IL2]{fontenc}
\usepackage{yfonts,fancyhdr}
\usepackage{a4wide}
\usepackage[english]{babel}
\usepackage{euscript}
\usepackage{amstext,amsbsy,amscd,amssymb}
\usepackage{amsmath,enumerate}
\usepackage{amsfonts}
\usepackage{mathrsfs,tensor,xy}
\usepackage{graphics}
\usepackage{graphicx}
\usepackage{microtype}
\usepackage{color}
\include{diagxy}

\let\rarr=\rightarrow

\let\veps=\varepsilon
\let\mcal=\mathcal
\let\mfrak=\mathfrak
\let\eus=\EuScript

\def\N{\mathbb{N}}
\def\Z{\mathbb{Z}}
\def\R{\mathbb{R}}
\def\C{\mathbb{C}}

\def\ad{\mathop {\rm ad} \nolimits}
\def\gr{\mathop {\rm gr} \nolimits}

\def\diag{\mathop {\rm diag} \nolimits}
\def\id{\mathop {\rm id} \nolimits}

\def\rank{\mathop {\rm rank} \nolimits}

\def\tr{\mathop {\rm tr} \nolimits}
\def\Res{\mathop {\rm Res} \nolimits}
\def\Ind{\mathop {\rm Ind} \nolimits}

\def\htt{\mathop {\rm ht} \nolimits}

\def\Fun{\mathop {\rm Fun} \nolimits}
\def\Vect{\mathop {\rm Vect} \nolimits}
\def\Pol{\mathop {\rm Pol} \nolimits}

\newsymbol\squares 1003

\long\def\proof #1{\noindent \emph{Proof.}\ #1 \hfill $\squares$
\medskip}

\newcounter{num}[section]
\numberwithin{equation}{section}
\numberwithin{num}{section}

\long\def\definition #1 {\refstepcounter{num} \noindent {\bf
Definition \thenum.} #1

\medskip}

\long\def\theorem #1{\refstepcounter{num} \noindent {\bf Theorem
\thenum.} #1

\medskip}

\long\def\lemma #1{\refstepcounter{num}  \noindent {\bf Lemma
\thenum.} #1

\medskip}

\long\def\proposition #1{\refstepcounter{num}  \noindent {\bf
Proposition \thenum.} #1

\medskip}

\long\def\corollary #1{\refstepcounter{num}  \noindent {\bf
Corollary \thenum.} #1

\medskip}

\makeatletter
\newcommand*\if@single[3]{%
  \setbox0\hbox{${\mathaccent"0362{#1}}^H$}%
  \setbox2\hbox{${\mathaccent"0362{\kern0pt#1}}^H$}%
  \ifdim\ht0=\ht2 #3\else #2\fi
  }
\newcommand*\rel@kern[1]{\kern#1\dimexpr\macc@kerna}
\newcommand*\widebar[1]{\@ifnextchar^{{\wide@bar{#1}{0}}}{\wide@bar{#1}{1}}}
\newcommand*\wide@bar[2]{\if@single{#1}{\wide@bar@{#1}{#2}{1}}{\wide@bar@{#1}{#2}{2}}}
\newcommand*\wide@bar@[3]{%
  \begingroup
  \def\mathaccent##1##2{%
    \if#32 \let\macc@nucleus\first@char \fi
    \setbox\z@\hbox{$\macc@style{\macc@nucleus}_{}$}%
    \setbox\tw@\hbox{$\macc@style{\macc@nucleus}{}_{}$}%
    \dimen@\wd\tw@
    \advance\dimen@-\wd\z@
    \divide\dimen@ 3
    \@tempdima\wd\tw@
    \advance\@tempdima-\scriptspace
    \divide\@tempdima 10
    \advance\dimen@-\@tempdima
    \ifdim\dimen@>\z@ \dimen@0pt\fi
    \rel@kern{0.6}\kern-\dimen@
    \if#31
      \overline{\rel@kern{-0.6}\kern\dimen@\macc@nucleus\rel@kern{0.4}\kern\dimen@}%
      \advance\dimen@0.4\dimexpr\macc@kerna
      \let\final@kern#2%
      \ifdim\dimen@<\z@ \let\final@kern1\fi
      \if\final@kern1 \kern-\dimen@\fi
    \else
      \overline{\rel@kern{-0.6}\kern\dimen@#1}%
    \fi
  }%
  \macc@depth\@ne
  \let\math@bgroup\@empty \let\math@egroup\macc@set@skewchar
  \mathsurround\z@ \frozen@everymath{\mathgroup\macc@group\relax}%
  \macc@set@skewchar\relax
  \let\mathaccentV\macc@nested@a
  \if#31
    \macc@nested@a\relax111{#1}%
  \else
    \def\gobble@till@marker##1\endmarker{}%
    \futurelet\first@char\gobble@till@marker#1\endmarker
    \ifcat\noexpand\first@char A\else
      \def\first@char{}%
    \fi
    \macc@nested@a\relax111{\first@char}%
  \fi
  \endgroup
}
\makeatother

\newcommand\rsmraise[1]{%
  \ifx#1\displaystyle .8\else
    \ifx#1\textstyle .8\else
      \ifx#1\scriptstyle .6\else
        .45%
      \fi
    \fi
  \fi}

\title{Geometric realizations of affine Kac-Moody algebras}

\author{Vyacheslav Futorny, Libor Křižka, Petr Somberg}

\AtEndDocument{\bigskip{\footnotesize%
  (V.\,Futorny) \textsc{Instituto de Matemática e Estatística, Universidade de Sa\~{o} Paulo, Caixa Postal 66281,\\ Sa\~{o} Paulo, CEP 05315-970, Brasil} \par
  \textit{E-mail address}: \texttt{futorny@ime.usp.br} \par
  \addvspace{\medskipamount}
  (L.\,Křižka) \textsc{Charles University, Faculty of Mathematics and Physics,
	Mathematical Institute, \\ Sokolovská 83, 180\,00 Praha 8, Czech Republic} \par
  \textit{E-mail address}: \texttt{krizka.libor@gmail.com} \par
  \addvspace{\medskipamount}
  (P.\,Somberg) \textsc{Charles University, Faculty of Mathematics and Physics,
	Mathematical Institute, \\ Sokolovská 83, 180\,00 Praha 8, Czech Republic} \par
  \textit{E-mail address}: \texttt{somberg@karlin.mff.cuni.cz} \par
}}

\date{}

\begin{document}

\maketitle

\begin{abstract}

The goal of the present paper is to obtain new free field realizations of affine Kac-Moody algebras
motivated by geometric representation theory for generalized flag manifolds of finite-dimensional
semisimple Lie groups. We provide an explicit construction of a large class of irreducible modules 
associated with certain parabolic subalgebras covering all known special cases.



\end{abstract}

\thispagestyle{empty}

\tableofcontents


\section*{Introduction}
\addcontentsline{toc}{section}{Introduction}

Classical free field realizations (or Fock space realizations) of affine Kac-Moody algebras were introduced by M.\,Wakimoto \cite{Wakimoto1986} for $\smash{\widehat{\mfrak{sl}}(2,\C)}$ and later generalized by B.\,Feigin and E.\,Frenkel \cite{Feigin-Frenkel1988} for an arbitrary affine Kac-Moody algebra $\widehat{\mfrak{g}}$ providing a construction of a family of $\widehat{\mfrak{g}}$-modules now called Wakimoto modules. These modules are generically irreducible and isomorphic to Verma modules induced from \emph{standard} Borel subalgebras of $\widehat{\mfrak{g}}$. On the other hand, generalized imaginary Verma modules  correspond to \emph{natural} Borel subalgebras (\cite{Jakobsen-Kac1985}, \cite{Futorny1994}, \cite{Futorny-Saifi1993}, \cite{Cox1994}). These modules have both finite-{} and infinite-dimensional weight spaces.  We refer the reader to \cite{Futorny1994}, \cite{Futorny1997} for their properties. 

A free field realization of a class of generalized imaginary Verma modules was obtained in \cite{Jakobsen-Kac1985} and \cite{Cox2005}. For any parabolic subalgebra $\mfrak{p}$ of a finite-dimensional simple Lie algebra $\mfrak{g}$ we define the \emph{natural} parabolic subalgebra $\mfrak{p}_{{\rm nat}}$ of $\widehat{\mfrak{g}}$ containing the \emph{natural} Borel subalgebra. Parabolic induction from a (continuous) representation $V$ of $\mfrak{p}_{{\rm nat}}$ leads to an induced $\widehat{\mfrak{g}}$-module called \emph{Generalized Imaginary Verma module}.
The particular classes of such modules are imaginary Verma modules, non-standard Verma type modules \cite{Futorny1997}, generalized
Wakimoto modules \cite{Cox-Futorny2004}, loop modules \cite{Chari1986} and generalized loop modules \cite{Bekkert-Benkart-Futorny-Kashuba2013}, \cite{Futorny-Kashuba2009}, \cite{Futorny-Kashuba2014}, \cite{Kashuba-Martins2014}.

The goal of our article is to give a uniform construction of geometrical origin for free field realizations of all Generalized Imaginary Verma modules for an arbitrary affine Kac-Moody algebras, thereby generalizing and unifying all known special constructions mentioned in the last paragraph.
Moreover, the results of \cite{Futorny-Kashuba2016}
imply that the Generalized Imaginary Verma modules are irreducible for 
a class (conjecturally for all) of inducing $\mfrak{p}_{{\rm nat}}$-modules
when the central charge is non-zero. 
Therefore, we have an explicit construction of a large new family of irreducible modules for affine Kac-Moody algebras, which should be considered as an analogue of the process well established for finite-dimensional semisimple Lie algebras and called localization of representations.

The structure of our article goes as follows. After basic preliminaries on affine
Kac-Moody algebras and their representation theory in Section \ref{sec1}, we discuss
in Section \ref{sec2} certain completions of infinite-dimensional Weyl algebras which
support the Generalized Imaginary Verma modules of affine Kac-Moody algebras.
In Section \ref{sec3} we present the general construction of Generalized Imaginary
Verma modules through their geometric realization. This is then
applied in Section \ref{sec4} to the affine Kac-Moody algebra $\smash{\widehat{\mfrak{sl}}(n+1,\C)}$
and its natural parabolic subalgebra $\mfrak{p}_{{\rm nat}}$ determined by the parabolic subalgebra $\mfrak{p} \simeq \C \oplus \mfrak{sl}(n,\C)$ of
$\mfrak{sl}(n+1,\C)$.

Our main results are Theorem \ref{thm:realization coordinates} which provides a geometric realization of affine Kac-Moody algebras, and Theorem \ref{thm:isomorphism} which shows an isomorphism between the geometric realization and the corresponding Generalized Imaginary Verma module. 

Following \cite{Frenkel-Ben-Zvi2004} and \cite{Frenkel2007} we work with topological versions of affine Kac-Moody algebras. Throughout the paper $\C$ is considered as a topological field with respect to the discrete topology. All topological vector spaces over $\C$ are supposed to be Hausdorff. We use the notation ${\mathbb Z}$ for the ring of integers, 
${\mathbb N}$ for natural numbers and ${\mathbb N}_0$ for natural numbers including zero.


\section{Affine Kac-Moody algebras}
\label{sec1}

Let $\mfrak{g}$ be a complex simple Lie algebra and let $\mfrak{h}$ be a Cartan subalgebra
of $\mfrak{g}$. We denote by $\Delta$ the root system of $\mfrak{g}$ with respect to $\mfrak{h}$,
by $\Delta^+$ the positive root system in $\Delta$ and by $\Pi \subset \Delta$ the set of
simple roots. We associate to the positive root system $\Delta^+$ the nilpotent Lie subalgebras
\begin{align}
  \mfrak{n} = \bigoplus_{\alpha \in \Delta^+} \mfrak{g}_\alpha \qquad \text{and} \qquad \widebar{\mfrak{n}} = \bigoplus_{\alpha \in \Delta^+} \mfrak{g}_{-\alpha}
\end{align}
and the solvable Lie subalgebras
\begin{align}
  \mfrak{b} = \mfrak{h} \oplus \mfrak{n} \qquad \text{and} \qquad \widebar{\mfrak{b}}= \mfrak{h} \oplus \widebar{\mfrak{n}}
\end{align}
of $\mfrak{g}$. The Lie algebras $\mfrak{b}$ and $\widebar{\mfrak{b}}$ are the standard and opposite standard Borel subalgebras of $\mfrak{g}$, respectively.

Let us consider a subset $\Sigma$ of $\Pi$ and denote by $\Delta_\Sigma$ the root subsystem in $\mfrak{h}^*$ generated by $\Sigma$. Then the standard parabolic subalgebra $\mfrak{p}$ and the opposite standard parabolic subalgebra $\widebar{\mfrak{p}}$ of $\mfrak{g}$ associated to $\Sigma$ are defined by
\begin{align}
  \mfrak{p} = \mfrak{l} \oplus \mfrak{u}  \qquad \text{and} \qquad \widebar{\mfrak{p}}= \mfrak{l} \oplus \widebar{\mfrak{u}},
\end{align}
where the reductive Levi subalgebra $\mfrak{l}$ of $\mfrak{p}$ and $\widebar{\mfrak{p}}$ is defined through
\begin{align}
  \mfrak{l}= \mfrak{h} \oplus \bigoplus_{\alpha \in \Delta_\Sigma} \mfrak{g}_\alpha
\end{align}
and the nilradical $\mfrak{u}$ of $\mfrak{p}$ and the opposite nilradical $\widebar{\mfrak{u}}$ of $\widebar{\mfrak{p}}$ are given by
\begin{align}
  \mfrak{u}= \bigoplus_{\alpha \in \Delta^+ \setminus \Delta^+_\Sigma} \mfrak{g}_\alpha \qquad \text{and} \qquad \widebar{\mfrak{u}}= \bigoplus_{\alpha \in \Delta^+ \setminus \Delta^+_\Sigma} \mfrak{g}_{-\alpha}.
\end{align}
We denote $\Delta^+ \setminus \Delta^+_\Sigma$ by $\Delta(\mfrak{u})$. Moreover, we have a triangular decomposition
\begin{align}
  \mfrak{g}= \widebar{\mfrak{u}} \oplus \mfrak{l} \oplus \mfrak{u}  \label{eq:triangular decomposition}
\end{align}
of the Lie algebra $\mfrak{g}$. Furthermore, we define the $\Sigma$-height $\htt_\Sigma(\alpha)$ of $\alpha \in \Delta$ by
\begin{align}
  \htt_\Sigma({\textstyle \sum_{i=1}^r} a_i\alpha_i)={\textstyle \sum^r_{i=1,\,\alpha_i \notin \Sigma}}\,  a_i,
\end{align}
where $r= \rank(\mfrak{g})$ and $\Pi=\{\alpha_1,\alpha_2,\dots,\alpha_r\}$. If we denote $k=\htt_\Sigma(\theta)$, where $\theta$ is the maximal root of $\mfrak{g}$ (by definition, the highest weight of its adjoint representation), then $\mfrak{g}$ is a $|k|$-graded Lie algebra with respect to the grading given by $\mfrak{g}_i= \bigoplus_{\alpha \in \Delta,\, \htt_\Sigma(\alpha)=i} \mfrak{g}_\alpha$ for $0 \neq i \in \Z$, and $\mfrak{g}_0 = \mfrak{h} \oplus \bigoplus_{\alpha \in \Delta,\, \htt_\Sigma(\alpha)=0} \mfrak{g}_\alpha$. Moreover, we have
\begin{align}
\widebar{\mfrak{u}}= \mfrak{g}_{-k} \oplus \dots \oplus \mfrak{g}_{-1}, \qquad \mfrak{l} =\mfrak{g}_0, \qquad  \mfrak{u}= \mfrak{g}_1 \oplus \dots \oplus \mfrak{g}_k.
\end{align}

Let $(\cdot\,,\cdot)\colon \mfrak{g} \otimes_\C \mfrak{g} \rarr \C$ be the $\mfrak{g}$-invariant symmetric bilinear form on $\mfrak{g}$ normalized by $(\theta,\theta)=2$, where $\theta$ is the maximal root of $\mfrak{g}$. Let us recall that we have
\begin{align}
  (\cdot\,,\cdot)= {1 \over 2h^\vee}\,(\cdot\,,\cdot)_\mfrak{g},
\end{align}
where $(\cdot\,,\cdot)_\mfrak{g} \colon \mfrak{g} \otimes_\C \mfrak{g} \rarr \C$ is the Cartan-Killing form on $\mfrak{g}$ and $h^\vee$ is the dual Coxeter number of $\mfrak{g}$.

The affine Kac-Moody algebra $\widehat{\mfrak{g}}$ associated to $\mfrak{g}$ is a universal central
extension of the formal loop algebra $\mfrak{g}(\!(t)\!)=\mfrak{g} \otimes_\C \C(\!(t)\!)$, i.e.\ $\widehat{\mfrak{g}}=\mfrak{g}(\!(t)\!)\oplus \C c$ with the commutation relations
\begin{align}
  [a\otimes f(t),b \otimes g(t)]=[a,b]\otimes f(t)g(t) - (a,b)\Res_{t=0}(f(t)dg(t))\,c, \label{eq:comm relations}
\end{align}
where $c$ is the central element of $\widehat{\mfrak{g}}$, $a,b \in \mfrak{g}$ and $f(t) ,g(t) \in \C(\!(t)\!)$. If we introduce the notation $a_n=a \otimes t^n$ for $a \in \mfrak{g}$ and $n \in \Z$, then \eqref{eq:comm relations} can be rewritten as
\begin{align}
  [a_m,b_n]=[a,b]_{m+n}+m(a,b)\delta_{m,-n}\,c \label{eq:comm relations modes}
\end{align}
for $m,n \in \Z$.
\medskip

We consider the natural Borel subalgebra $\smash{\widehat{\mfrak{b}}_{\rm nat}}$ of $\widehat{\mfrak{g}}$ defined by
\begin{align}
  \widehat{\mfrak{b}}_{\rm nat} = \widehat{\mfrak{h}}_{\rm nat} \oplus \widehat{\mfrak{n}}_{\rm nat},
\end{align}
where the Cartan subalgebra $\widehat{\mfrak{h}}_{\rm nat}$ is given by
\begin{align}
  \widehat{\mfrak{h}}_{\rm nat} = \mfrak{h} \otimes_\C \C 1 \oplus \C c
\end{align}
and the nilradical $\widehat{\mfrak{n}}_{\rm nat}$ of $\widehat{\mfrak{b}}_{\rm nat}$ and the opposite nilradical $\widehat{\widebar{\mfrak{n}}}_{\rm nat}$ are
\begin{align}
  \widehat{\mfrak{n}}_{\rm nat} = \mfrak{n} \otimes_\C \C(\!(t)\!) \oplus \mfrak{h} \otimes_\C t\C[[t]] \qquad \text{and} \qquad \widehat{\widebar{\mfrak{n}}}_{\rm nat}= \widebar{\mfrak{n}} \otimes_\C \C(\!(t)\!) \oplus \mfrak{h} \otimes_\C t^{-1}\C[t^{-1}].
\end{align}
Moreover, we have a triangular decomposition
\begin{align}
  \widehat{\mfrak{g}}= \widehat{\widebar{\mfrak{n}}}_{\rm nat} \oplus \widehat{\mfrak{h}}_{\rm nat} \oplus \widehat{\mfrak{n}}_{\rm nat}
\end{align}
of the Lie algebra $\widehat{\mfrak{g}}$.
In addition, we introduce the natural parabolic subalgebra $\mfrak{p}_{\rm nat}$ of $\widehat{\mfrak{g}}$ associated to a parabolic subalgebra $\mfrak{p}$ of $\mfrak{g}$ through
\begin{align}
  \mfrak{p}_{\rm nat} = \mfrak{l}_{\rm nat} \oplus \mfrak{u}_{\rm nat},
\end{align}
where the reductive Levi subalgebra $\mfrak{l}_{\rm nat}$ of $\mfrak{p}_{\rm nat}$ is defined by
\begin{align}
\mfrak{l}_{\rm nat}=\mfrak{l}\otimes_\C \C(\!(t)\!) \oplus \C c
\end{align}
and the nilradical $\mfrak{u}_{\rm nat}$ of $\mfrak{p}_{\rm nat}$ and the opposite nilradical $\widebar{\mfrak{u}}_{\rm nat}$ are given by
\begin{align}
  \mfrak{u}_{\rm nat} = \mfrak{u}\otimes_\C \C(\!(t)\!) \qquad \text{and} \qquad \widebar{\mfrak{u}}_{\rm nat}= \widebar{\mfrak{u}}\otimes_\C \C(\!(t)\!).
\end{align}
Therefore, we have a triangular decomposition
\begin{align}
  \widehat{\mfrak{g}}= \widebar{\mfrak{u}}_{\rm nat} \oplus  \mfrak{l}_{\rm nat} \oplus \mfrak{u}_{\rm nat}
\end{align}
of the Lie algebra $\widehat{\mfrak{g}}$.
\medskip

\definition{Let $\sigma \colon \mfrak{p}_{\rm nat} \rarr \mfrak{gl}(V)$ be a $\mfrak{p}_{\rm nat}$-module such that $\sigma(c)=k\cdot \id_V$ for $k\in \C$. Then the generalized imaginary Verma module at level $k$ is the induced module
\begin{align}
  \mathbb{M}_{\sigma,k,\mfrak{p}}(V) = \Ind^{\widehat{\mfrak{g}}}_{\mfrak{p}_{\rm nat}}\!V \equiv U(\widehat{\mfrak{g}}) \otimes_{U(\mfrak{p}_{\rm nat})}\!V \simeq U(\widebar{\mfrak{u}}_{\rm nat}) \otimes_\C\! V,
\end{align}
where the last isomorphism of vector spaces follows from Poincar\'e–Birkhoff–Witt theorem.}
\vspace{-2mm}


\section{Weyl algebras and their completions}
\label{sec2}

In this section we introduce a formalism for infinite-dimensional Weyl algebras and
define several of their completions. This will enable us to construct generalized imaginary Verma modules for affine Kac-Moody algebras.


\subsection{Weyl algebras in the infinite-dimensional setting}

Let us consider the commutative $\C$-algebra $\mcal{K}=\C(\!(t)\!)$ with the $\C$-subalgebra $\mcal{O}=\C[[t]]$. Let $\Omega_\mcal{K}=\C(\!(t)\!)\,dt$ and $\Omega_\mcal{O}=\C[[t]]\,dt$ be the modules of Kähler differentials. For a finite-dimensional complex vector space $V$ we define the infinite-dimensional complex vector spaces $\mcal{K}(V)=V\otimes_\C \mcal{K}$ and $\Omega_\mcal{K}(V^*)=V^*\otimes_\C \Omega_\mcal{K}$. The pairing $(\cdot\,,\cdot) \colon \Omega_\mcal{K}(V^*) \otimes_\C \mcal{K}(V) \rarr \C$ defined by
\begin{align}
  (\alpha \otimes f(t)dt,v \otimes g(t))=\alpha(v) \Res_{t=0}(g(t)f(t)dt), \label{eq:non-degenerate form}
\end{align}
where $\alpha \in V^*$, $v \in V$ and $f(t),g(t) \in \mcal{K}$, allows to identify the restricted dual space to $\mcal{K}(V)$ with the vector space $\Omega_\mcal{K}(V^*)$, and vice versa. Moreover, the pairing \eqref{eq:non-degenerate form} gives us a skew-symmetric non-degenerate bilinear form $\langle \cdot\,,\cdot \rangle$ on $\Omega_\mcal{K}(V^*) \oplus \mcal{K}(V)$ defined by
\begin{align}
\langle \alpha \otimes f(t)dt,v \otimes g(t) \rangle = - \langle v \otimes g(t), \alpha \otimes f(t)dt \rangle = \alpha(v) \Res_{t=0} (g(t)f(t)dt)
\end{align}
for $\alpha \in V^*$, $v \in V$, $f(t), g(t) \in \mcal{K}$, and
\begin{align}
  \langle v \otimes f(t), w \otimes g(t) \rangle = \langle \alpha \otimes f(t)dt, \beta \otimes g(t)dt \rangle = 0
\end{align}
for $\alpha,\beta \in V^*$, $v,w \in V$, $f(t),g(t) \in \mcal{K}$. Then the Weyl algebras $\eus{A}_{\mcal{K}(V)}$ and $\eus{A}_{\Omega_\mcal{K}(V^*)}$ are given by
\begin{align}
  \eus{A}_{\mcal{K}(V)} = T(\Omega_\mcal{K}(V^*) \oplus \mcal{K}(V))/I_{\mcal{K}(V)} \qquad \text{and} \qquad \eus{A}_{\Omega_\mcal{K}(V^*)} = T(\Omega_\mcal{K}(V^*) \oplus \mcal{K}(V))/I_{\Omega_\mcal{K}(V^*)},
\end{align}
where $I_{\mcal{K}(V)}$ and $I_{\Omega_\mcal{K}(V^*)}$ denote the two-sided ideals of the tensor algebra $T(\Omega_\mcal{K}(V^*) \oplus \mcal{K}(V))$ generated by $a \otimes b - b \otimes a +\langle a,b \rangle \cdot 1$ for all $a \in \Omega_\mcal{K}(V^*)$, $b \in \mcal{K}(V)$ and by $a \otimes b - b \otimes a -\langle a,b \rangle \cdot 1$ for all $a \in \Omega_\mcal{K}(V^*)$, $b \in \mcal{K}(V)$, respectively. The $\C$-algebras of polynomials on $\mcal{K}(V)$ and $\Omega_\mcal{K}(V^*)$ we define by $\Pol \mcal{K}(V) = S(\Omega_\mcal{K}(V^*))$ and $\Pol \Omega_\mcal{K}(V^*) = S(\mcal{K}(V))$, respectively.

Let $\mcal{L}$ and $\mcal{L}^{\rm c}$ be complementary Lagrangian (maximal isotropic) subspaces of $\Omega_\mcal{K}(V^*) \oplus \mcal{K}(V)$, i.e., we have $\Omega_\mcal{K}(V^*) \oplus \mcal{K}(V) = \mcal{L} \oplus \mcal{L}^{\rm c}$. The most interesting case for us will be $\mcal{L} = \Omega_\mcal{K}(V^*)$ and $\mcal{L}^{\rm c} = \mcal{K}(V)$.

The symmetric algebra $S(\mcal{L})$ is a subalgebra of $\eus{A}_{\mcal{K}(V)}$, since the elements of $\mcal{L}$ commute in $\eus{A}_{\mcal{K}(V)}$. In fact, it is a maximal commutative $\C$-subalgebra of $\eus{A}_{\mcal{K}(V)}$. We consider the induced $\eus{A}_{\mcal{K}(V)}$-module
\begin{align}
  \Ind_{S(\mcal{L})}^{\eus{A}_{\mcal{K}(V)}} \C \simeq S(\mcal{L}^{\rm c}), \label{eq:induced module}
\end{align}
where $\C$ is the trivial $S(\mcal{L})$-module. It follows from \eqref{eq:induced module} that the
induced $\eus{A}_{\mcal{K}(V)}$-module has a structure of a commutative $\C$-algebra. We denote by
$M_\mcal{L}$ the commutative $\C$-algebra which is the completion of $S(\mcal{L}^{\rm c})$ with
respect to the linear topology on $S(\mcal{L}^{\rm c})$ in which the basis of open neighborhoods of $0$ are the subspaces $\mcal{I}_n$ for $n\in \Z$, where $\mcal{I}_n$ is the ideal of $S(\mcal{L}^{\rm c})$ generated by $\mcal{L}^{\rm c}\cap (V^*\otimes_\C t^n\Omega_\mcal{O})$. We can extend the action of the Weyl algebra $\eus{A}_{\mcal{K}(V)}$ to $M_\mcal{L}$.
\medskip

Our next step is to pass to a completion of the Weyl algebra $\smash{\eus{A}_{\mcal{K}(V)}}$,
because $\smash{\eus{A}_{\mcal{K}(V)}}$ is not sufficiently large for our considerations.
Let us denote by $\Fun\mcal{K}(V)$ the completion of the commutative $\C$-algebra $\Pol \mcal{K}(V)$ with respect to the linear topology on $\Pol \mcal{K}(V)$ in which the basis of open neighborhoods of $0$ are the subspaces $\mcal{J}_n$ for $n \in \Z$, where $\mcal{J}_n$ is the ideal of $\Pol \mcal{K}(V)$ generated by $V^* \otimes_\C t^n\Omega_\mcal{O}$.
Consequently, we have $\Fun \mcal{K}(V)= M_{\mcal{K}(V)}$.
Then a vector field on $\mcal{K}(V)$ is by definition a continuous $\C$-linear endomorphism $\xi$ of $\Fun \mcal{K}(V)$ which satisfies the Leibniz rule
\begin{align}
  \xi(fg)=\xi(f)g+f\xi(g)
\end{align}
for all $f,g \in \Fun \mcal{K}(V)$. The vector space of all vector fields is naturally a topological Lie algebra, which we denote by $\Vect \mcal{K}(V)$.
Moreover, there is a split short exact sequence
\begin{align}
  0 \rarr \Fun \mcal{K}(V) \rarr \Vect \mcal{K}(V) \oplus \Fun \mcal{K}(V) \rarr \Vect \mcal{K}(V) \rarr 0
\end{align}
of topological Lie algebras. We define the completed Weyl algebra $\smash{\eus{A}^\sharp_{\mcal{K}(V)}}$ as the associative $\C$-algebra generated by a subalgebra $i \colon \Fun\mcal{K}(V) \rarr \smash{\eus{A}^\sharp_{\mcal{K}(V)}}$ and a Lie subalgebra $j \colon  \Vect\mcal{K}(V) \oplus \Fun \mcal{K}(V) \rarr \smash{\eus{A}^\sharp_{\mcal{K}(V)}}$, with the relation
\begin{align}
  [j(\xi+f),i(g)]=i(\xi(g))
\end{align}
for all $f, g \in \Fun \mcal{K}(V)$ and $\xi \in \Vect \mcal{K}(V)$. We obtain that $\Fun \mcal{K}(V)$ is an $\smash{\eus{A}^\sharp_{\mcal{K}(V)}}$-module.
\medskip

Let $\{f_\alpha;\, \alpha \in \Delta(\mfrak{u})\}$ be a basis of $\widebar{\mfrak{u}}$. Further, let $\{x_\alpha;\, \alpha \in \Delta(\mfrak{u})\}$ be the linear coordinate functions on $\widebar{\mfrak{u}}$ with respect to the given basis of $\widebar{\mfrak{u}}$. Then the set $\{f_\alpha \otimes t^n;\, \alpha \in \Delta(\mfrak{u}),\, n \in \Z\}$ forms a topological basis of $\mcal{K}(\widebar{\mfrak{u}})=\widebar{\mfrak{u}}_{{\rm nat}}$, and the set $\{x_\alpha \otimes t^{-n-1}dt;\, \alpha \in \Delta(\mfrak{u}),\, n \in \Z\}$ forms a dual topological basis of $\Omega_\mcal{K}(\widebar{\mfrak{u}}^*) \simeq (\widebar{\mfrak{u}}_{{\rm nat}})^*$ with respect to the pairing \eqref{eq:non-degenerate form}, i.e.\ we have
\begin{align}
  (x_\alpha \otimes t^{-n-1}dt,f_\beta \otimes t^m)=x_\alpha(f_\beta)\Res_{t=0} t^{m-n-1}dt= \delta_{\alpha\beta} \delta_{nm}
\end{align}
for all $\alpha,\beta \in \Delta(\mfrak{u})$ and $m,n \in \Z$. If we denote $x_{\alpha,n}=x_\alpha \otimes t^{-n-1}dt$ and $\partial_{x_{\alpha,n}}=f_\alpha \otimes t^n$ for $\alpha \in \Delta(\mfrak{u})$ and $n \in \Z$, then the two-sided ideal $I_{\mcal{K}(\widebar{\mfrak{u}})}$ is generated by elements
\begin{align}
\bigg(\sum_{n \in \Z} a_nx_{\alpha,n}\bigg)\! \otimes \!\bigg(\sum_{m\in \Z} b_m\partial_{x_{\beta,m}}\! \bigg) - \bigg(\sum_{m\in \Z} b_m\partial_{x_{\beta,m}}\! \bigg) \!\otimes\! \bigg(\sum_{n \in \Z} a_nx_{\alpha,n}\bigg)+\delta_{\alpha\beta}\bigg(\sum_{n\in \Z} a_nb_n\bigg)\!\cdot 1
\end{align}
and coincides with the canonical commutation relations
\begin{align}
[x_{\alpha,n},\partial_{x_{\beta,m}}]=-\delta_{\alpha\beta}\delta_{nm}
\end{align}
for all $\alpha,\beta \in \Delta(\mfrak{u})$ and $m,n \in \Z$. Therefore, we obtain that the Weyl algebra $\eus{A}_{\mcal{K}(\widebar{\mfrak{u}})}$ is topologically generated by $\{x_{\alpha,n},\partial_{x_{\alpha,n}};\, \alpha \in \Delta(\mfrak{u}),\, n \in \Z\}$ with the canonical commutation relations. We have $\smash{M_{\Omega_\mcal{K}(\widebar{\mfrak{u}}^*)}}=\Pol \Omega_\mcal{K}(\widebar{\mfrak{u}}^*)$ with the discrete topology, and $M_{\mcal{K}(\widebar{\mfrak{u}})}=\Fun \mcal{K}(\widebar{\mfrak{u}})$ with the linear topology in which the basis of open neighbourhoods of $0$ are the subspaces $\mcal{J}_n$ for $n \in \Z$, where $\mcal{J}_n$ is the ideal of $\Fun \mcal{K}(\widebar{\mfrak{u}})$ generated by $\widebar{\mfrak{u}}^* \otimes_\C t^n\Omega_\mcal{O}$.

For later purposes, the completed Weyl algebra $\smash{\eus{A}^\sharp_{\mcal{K}(\widebar{\mfrak{u}})}}$  will be denoted by $\eus{A}^{\mfrak{g},\mfrak{p}}$.


\subsection{The local extension}

For our purposes we may replace the completed Weyl algebra $\eus{A}^{\mfrak{g},\mfrak{p}}$, which is a very large topological algebra, by a relatively small local part. So let us consider the vector space $\eus{F}^{\mfrak{g},\mfrak{p}}_{{\rm loc}}$ of local functions on $\mcal{K}(\widebar{\mfrak{u}})$ spanned by the Fourier coefficients of the form
\begin{align}
  \Res_{z=0} P(a_\alpha^*(z),\partial_z a_\alpha^*(z),\dots)f(z)dz,
\end{align}
where $P(a_\alpha^*(z),\partial_z a_\alpha^*(z),\dots)$ is a differential polynomial in $a^*_\alpha(z)$ for $\alpha \in \Delta(\mfrak{u})$ and $f(z) \in \C(\!(z)\!)$. Further, let $\eus{T}^{\mfrak{g},\mfrak{p}}_{{\rm loc}}$ be the vector space of local vector fields on $\mcal{K}(\widebar{\mfrak{u}})$ spanned by the Fourier coefficients of the form
\begin{align}
  \Res_{z=0} P(a_\alpha^*(z),\partial_z a_\alpha^*(z),\dots)a_\beta(z)f(z)dz, \label{eq:vector fields local}
\end{align}
where we used the fact that the Fourier coefficients of the form
\begin{align}
  \Res_{z=0} P(a_\alpha^*(z),\partial_z a_\alpha^*(z),\dots)\partial_z^m a_\beta(z)f(z)dz
\end{align}
for $m>0$ may be expressed as linear combinations of \eqref{eq:vector fields local}. Finally, let us consider the vector space $\smash{\eus{A}^{\mfrak{g},\mfrak{p}}_{\leq m,{\rm loc}}}$ of local differential operators of order at most $m \in \N_0$ on $\mcal{K}(\widebar{\mfrak{u}})$ spanned by the Fourier coefficients of the form
\begin{align}
  \Res_{z=0} P(a_\alpha^*(z),\partial_z a_\alpha^*(z),\dots)Q(a_\beta(z),\partial_z a_\beta(z),\dots)f(z)dz,
\end{align}
where $Q(a_\beta(z),\partial_z a_\beta(z),\dots)$ is a differential polynomial in $a_\beta(z)$ for $\beta \in \Delta(\mfrak{u})$ of degree at most $m$. Then $\smash{\eus{A}^{\mfrak{g},\mfrak{p}}_{\leq m,{\rm loc}}}$ for $m \in \N_0$ are topological Lie algebras and we get a short exact sequence
\begin{align}
  0 \rarr \eus{F}^{\mfrak{g},\mfrak{p}}_{{\rm loc}} \rarr \eus{A}^{\mfrak{g},\mfrak{p}}_{\leq 1,{\rm loc}} \rarr \eus{T}^{\mfrak{g},\mfrak{p}}_{{\rm loc}} \rarr 0 \label{eq:exact sequence}
\end{align}
of topological Lie algebras, which has a canonical splitting. Furthermore, we define a topological Lie algebra $\eus{A}^{\mfrak{g},\mfrak{p}}_{{\rm loc}}$ by
\begin{align}
  \eus{A}^{\mfrak{g},\mfrak{p}}_{{\rm loc}} = \bigcup_{m \in \N_0} \eus{A}^{\mfrak{g},\mfrak{p}}_{\leq m, {\rm loc}}.
\end{align}
Let us note that $\eus{A}^{\mfrak{g},\mfrak{p}}_{{\rm loc}}$ is only a Lie algebra and not an associative algebra. However, there is a construction which associates to $\eus{A}^{\mfrak{g},\mfrak{p}}_{{\rm loc}}$ the completed Weyl algebra $\eus{A}^{\mfrak{g},\mfrak{p}}$.

Now, from the construction we have that $M_{\mcal{K}(\widebar{\mfrak{u}})}$ is an $\eus{A}^{\mfrak{g},\mfrak{p}}$-module, however $\smash{M_{\Omega_\mcal{K}(\widebar{\mfrak{u}}^*)}}$ is not an $\eus{A}^{\mfrak{g},\mfrak{p}}$-module. Therefore, we need to consider a different completion of the Weyl algebra $\eus{A}_{\mcal{K}(\widebar{\mfrak{u}})}$. Let us note that there are more completions which are closely related to different normal orderings of differential operators. We restrict our attention to a completion suitable for a free field realization of generalized imaginary Verma modules.

Let us consider the vector space $\smash{\eus{A}^{\mfrak{g},\mfrak{p},{\rm op}}_{\leq m,{\rm loc}}}$ of local differential operators of order at most $m \in \N_0$ on $\mcal{K}(\widebar{\mfrak{u}})$ spanned by the Fourier coefficients of the form
\begin{align}
  \Res_{z=0} Q(a_\beta(z),\partial_z a_\beta(z),\dots)P(a_\alpha^*(z),\partial_z a_\alpha^*(z),\dots)f(z)dz,
\end{align}
where $Q(a_\beta(z),\partial_z a_\beta(z),\dots)$ is a differential polynomial in $a_\beta(z)$ for $\beta \in \Delta(\mfrak{u})$ of degree at most $m$. Then $\smash{\eus{A}^{\mfrak{g},\mfrak{p},{\rm op}}_{\leq m,{\rm loc}}}$ for $m\in \N_0$ are topological Lie algebras. To construct a short exact sequence of topological Lie algebras similar to \eqref{eq:exact sequence}, we define a mapping $\varphi \colon \smash{\eus{A}^{\mfrak{g},\mfrak{p},{\rm op}}_{\leq 1,{\rm loc}}} \rarr \smash{\eus{A}^{\mfrak{g},\mfrak{p}}_{\leq 1,{\rm loc}}}$ by
\begin{align}
  \varphi(\Res_{z=0} Q(a_\beta(z),\dots)P(a_\alpha^*(z),\dots)f(z)dz)= \Res_{z=0} P(a_\alpha^*(z),\dots)Q(a_\beta(z),\dots)f(z)dz.
\end{align}
Then we have
\begin{multline*}
  \varphi([\Res_{z=0} Q_1(z)P_1(z)f_1(z)dz,\Res_{w=0} Q_2(w)P_2(w)f_2(w)dw])\\ =\varphi\big(\!\Res_{z=0,w=0}\big(Q_1(z)[P_1(z),Q_2(w)]P_2(w) +Q_2(w)[Q_1(z),P_2(w)]P_1(z)\big)f_1(z)f_2(w)dzdw\big) \\
  =\Res_{z=0,w=0}\big(P_2(w)[P_1(z),Q_2(w)]Q_1(z)
   +P_1(z)[Q_1(z),P_2(w)]Q_2(w)\big)f_1(z)f_2(w)dzdw \\
  = [\Res_{z=0} P_1(z)Q_1(z)f_1(z)dz,\Res_{w=0} P_2(w)Q_2(w)f_2(w)dw]\\
  =[\varphi(\Res_{z=0} Q_1(z)P_1(z)f_1(z)dz),\varphi(\Res_{w=0} Q_2(w)P_2(w)f_2(w)dw)]),
\end{multline*}
where $Q_1(z),Q_2(z)$ are differential polynomials in $a_\beta(z)$ for $\beta \in \Delta(\mfrak{u})$ of degree at most 1 and $P_1(z),P_2(z)$ are differential polynomials in $a^*_\alpha(z)$ for $\alpha \in \Delta(\mfrak{u})$,
which implies that the mapping $\varphi \colon \smash{\eus{A}^{\mfrak{g},\mfrak{p},{\rm op}}_{\leq 1,{\rm loc}}} \rarr \smash{\eus{A}^{\mfrak{g},\mfrak{p}}_{\leq 1,{\rm loc}}}$ is an isomorphism of topological Lie algebras. Hence, we get a short exact sequence
\begin{align}
  0 \rarr \eus{F}^{\mfrak{g},\mfrak{p}}_{{\rm loc}} \rarr \eus{A}^{\mfrak{g},\mfrak{p},{\rm op}}_{\leq 1,{\rm loc}} \rarr \eus{T}^{\mfrak{g},\mfrak{p}}_{{\rm loc}} \rarr 0
\end{align}
of topological Lie algebras, which has a canonical splitting. Furthermore, we define a topological Lie algebra $\eus{A}^{\mfrak{g},\mfrak{p},{\rm op}}_{{\rm loc}}$ by
\begin{align}
  \eus{A}^{\mfrak{g},\mfrak{p},{\rm op}}_{{\rm loc}} = \bigcup_{m \in \N_0} \eus{A}^{\mfrak{g},\mfrak{p},{\rm op}}_{\leq m, {\rm loc}}.
\end{align}
Let us note that $\eus{A}^{\mfrak{g},\mfrak{p},{\rm op}}_{{\rm loc}}$ is again only a Lie algebra and not an associative algebra. However, there is a construction which associates to $\eus{A}^{\mfrak{g},\mfrak{p},{\rm op}}_{{\rm loc}}$ the completed Weyl algebra $\eus{A}^{\mfrak{g},\mfrak{p},{\rm op}}$.

Furthermore, let us consider the vector space $\eus{J}^{\mfrak{g},\mfrak{p}}_{{\rm loc}}\subset \eus{F}^{\mfrak{g},\mfrak{p}}\,\widehat{\otimes}\, \mfrak{p}_{{\rm nat}}$ spanned by the Fourier coefficients of the form
\begin{align}
  \Res_{z=0}P(a_\alpha^*(z),\partial_z a_\alpha^*(z),\dots)h(z)f(z)dz
\end{align}
and
\begin{align}
  \Res_{z=0}P(a_\alpha^*(z),\partial_z a_\alpha^*(z),\dots)cf(z)dz,
\end{align}
where $P(a_\alpha^*(z),\partial_z a_\alpha^*(z),\dots)$ is a differential polynomial in  $a^*_\alpha(z)$ for $\alpha \in \Delta(\mfrak{u})$, $h \in \mfrak{p}$ and $f(z) \in \C(\!(z)\!)$. Moreover, we have that $\eus{J}^{\mfrak{g},\mfrak{p}}_{{\rm loc}}$ is a topological Lie algebra and induces a natural structure of a topological Lie algebra on a semidirect sum
\begin{align}
  \eus{R}^{\mfrak{g},\mfrak{p}}_{{\rm loc}} = \smash{\eus{A}^{\mfrak{g},\mfrak{p},{\rm op}}_{\leq 1,{\rm loc}}} \oplus \eus{J}^{\mfrak{g},\mfrak{p}}_{{\rm loc}}.
\end{align}
This topological Lie algebra plays the key role in a free field realization of generalized imaginary Verma modules as we shall observe in the next section.


\section{Generalized Imaginary Verma modules}
\label{sec3}

In this section we will give a construction of a new family of free field realizations 
of affine Kac-Moody algebras, based on completed infinite-dimensional Weyl algebra and 
the generalized imaginary Verma modules.


\subsection{A geometric realization of affine Kac-Moody algebras}

Let $R$ be an algebra over $\C$. Then an $R$-valued formal power series (or formal distribution) in the variables $z_1,z_2,\dots,z_n$ is a series
\begin{align}
  a(z_1,\dots,z_n) = \sum_{m_1,\dots,m_n \in \Z} a_{m_1 \dots m_n} z_1^{m_1} \dots z_n^{m_n},
\end{align}
where $a_{m_1 \dots m_n}  \in R$. The vector space of all formal power series is denoted by $R[[z_1^{\pm 1},\dots,z_n^{\pm 1}]]$.

An important example of a $\C$-valued formal power series in two variables $z,w$ is the formal delta function $\delta(z-w)$ defined by
\begin{align}
  \delta(z-w) = \sum_{m \in \Z} z^m w^{-m-1}.  \label{eq:delta function}
\end{align}
The useful properties of the formal delta function are summarized in the following proposition which is standard (cf.\ \cite{Frenkel-Ben-Zvi2004}).
\medskip

\proposition{\label{prop:delta function}
We have
\begin{enumerate}
  \item[1)] $\delta(z-w)=\delta(w-z)$,
  \item[2)] $\partial_z\delta(z-w)= -\partial_w \delta(z-w)$,
  \item[3)] $a(z)\delta(z-w)=a(w)\delta(z-w)$,
  \item[4)] $a(z)\partial_w \delta(z-w) = a(w)\partial_w \delta(z-w)+ (\partial_w a(w))\delta(z-w)$,
  \item[5)] $(z-w)^{n+1}\partial_w^n \delta(z-w)=0$,
  \item[6)] $\Res_{z=0} a(z)\delta(z-w)=a(w)$
\end{enumerate}
for all $n \in \N_0$ and $a(z) \in R[[z^{\pm 1}]]$.}








We will use the notation of the previous sections. Let $\{f_\alpha;\, \alpha \in \Delta(\mfrak{u})\}$ be a basis of $\widebar{\mfrak{u}}$ and let $\{x_\alpha;\, \alpha \in \Delta(\mfrak{u})\}$ be the linear coordinate functions on $\widebar{\mfrak{u}}$ with respect to the given basis of $\widebar{\mfrak{u}}$.
For $a \in \mfrak{g}$ we define the formal distribution $a(z) \in \widehat{\mfrak{g}}[[z^{\pm 1}]]$ by
\begin{align}
  a(z)= \sum_{n\in \Z} a_n z^{-n-1}, \label{eq:formal dist algebra}
\end{align}
where $a_n = a \otimes t^n$ for $n \in \Z$. Then the commutation relations \eqref{eq:comm relations modes} can be equivalently written as
\begin{align}
  [a(z),b(w)]= [a,b](w) \delta(z-w) + (a,b)c\,\partial_w \delta(z-w) \label{eq:comm relations series}
\end{align}
for all $a,b \in \mfrak{g}$. Further, we also introduce the formal distributions $a_\alpha(z), a_\alpha^*(z) \in \eus{A}_{\mcal{K}(\widebar{\mfrak{u}})}[[z^{\pm 1}]]$ by
\begin{align}
  a_\alpha(z)= \sum_{n\in \Z} a_{\alpha,n}z^{-n-1} \qquad \text{and} \qquad a_\alpha^*(z) = \sum_{n\in \Z} a^*_{\alpha,n}z^{-n}, \label{eq:formal dist a,a*}
\end{align}
where $a_{\alpha,n}=\partial_{x_{\alpha,n}}$ and $a^*_{\alpha,n}=x_{\alpha,-n}$, for $\alpha \in \Delta(\mfrak{u})$.
\medskip

\proposition{We have
\begin{align}
  [a_\alpha(z),a_\beta(w)]=0, \qquad [a_\alpha(z),a^*_\beta(w)]=\delta_{\alpha\beta}\delta(z-w), \qquad [a^*_\alpha(z),a^*_\beta(w)]=0
\end{align}
for $\alpha,\beta \in \Delta(\mfrak{u})$.}

\proof{We have
\begin{align*}
  [a_\alpha(z),a^*_\beta(w)]&=\sum_{n\in \Z} \sum_{m\in \Z}[\partial_{x_{\alpha,n}},x_{\beta,m}]z^{-n-1}w^m = \sum_{m\in \Z} \sum_{n\in \Z} \delta_{\alpha\beta} \delta_{mn}z^{-n-1}w^m \\
  &=\delta_{\alpha\beta} \sum_{m\in \Z} z^{-m-1}w^m = \delta_{\alpha\beta} \delta(z-w)
\end{align*}
for $\alpha,\beta \in \Delta(\mfrak{u})$. The other two commutation relations have analogous proof.}

The following theorem is our first main result. It gives an embedding of the affine Kac-Moody algebra $\widehat{\mfrak{g}}$ into the topological Lie algebra $\eus{R}^{\mfrak{g},\mfrak{p}}_{{\rm loc}}$.
\medskip

\theorem{\label{thm:realization coordinates}
The embedding of $\widehat{\mfrak{g}}$ into $\eus{R}^{\mfrak{g},\mfrak{p}}_{{\rm loc}}$ is given by
\begin{align}
\pi(a)= -\sum_{\alpha \in \Delta(\mfrak{u})} \sum_{n \in \Z} \partial_{x_{\alpha,n}} \bigg[{\ad(u(x))e^{\ad(u(x))} \over e^{\ad(u(x))}-{\rm id}_{\widebar{\mfrak{u}}_{{\rm nat}}}}\,(e^{-\ad(u(x))}a)_{\widebar{\mfrak{u}}_{{\rm nat}}}\bigg]_{\alpha,n} + (e^{-\ad(u(x))}a)_{\mfrak{p}_{\rm nat}}  \label{eq:explicit realization coordinates}
\end{align}
for all $a\in \widehat{\mfrak{g}}$, where
\begin{align}
  u(x)=\sum_{\alpha \in \Delta(\mfrak{u})} \sum_{n \in \Z} x_{\alpha,n} f_{\alpha,n}
\end{align}
and $[b]_{\alpha,n}$ is the $(\alpha,n)$-th coordinate of $b \in \widebar{\mfrak{u}}_{{\rm nat}}$ with respect to the topological basis $\{f_{\alpha,n};\, \alpha \in \Delta(\mfrak{u}),\, n \in \Z\}$ of $\widebar{\mfrak{u}}_{{\rm nat}}$.}

\proof{It follows from the equivalent reformulation in Theorem \ref{thm:realization power series}.}

We will rewrite the formula \eqref{eq:explicit realization coordinates} using the formal distributions \eqref{eq:formal dist algebra} and \eqref{eq:formal dist a,a*} into an equivalent and  more compact form.
Let us denote by $\eus{P}^{\mfrak{g},\mfrak{p}}_{{\rm loc}}(z)$ the vector space of all polynomials in $a^*_\alpha(z)$ for $\alpha \in \Delta(\mfrak{u})$, by $\eus{F}^{\mfrak{g},\mfrak{p}}_{{\rm loc}}(z)$ the vector space of all differential polynomials in $a^*_\alpha(z)$ for $\alpha \in \Delta(\mfrak{u})$, and by $\eus{C}^\mfrak{g}_{{\rm loc}}(z)$ the vector space of all formal distributions of the form $a(z)$ for $a \in \mfrak{g}$.
We define a formal power series $u(z) \in \widebar{\mfrak{u}} \otimes_\C \eus{P}^{\mfrak{g},\mfrak{p}}_{{\rm loc}}(z)$ by
\begin{align}
  u(z) = \sum_{\alpha \in \Delta(\mfrak{u})} a_\alpha^*(z)f_\alpha. \label{eq:u(z) def}
\end{align}
Since the mapping $a \mapsto a(z)$ from $\mfrak{g}$ to $\eus{C}^\mfrak{g}_{{\rm loc}}(z)$ is an isomorphism of vector spaces, we define a structure of a Lie algebra on $\eus{C}^\mfrak{g}_{{\rm loc}}(z)$ via this isomorphism. Therefore, we obtain that $\mfrak{g} \otimes_\C \eus{P}^{\mfrak{g},\mfrak{p}}_{{\rm loc}}(z)$ and $\eus{C}^\mfrak{g}_{{\rm loc}}(z) \otimes_\C \eus{P}^{\mfrak{g},\mfrak{p}}_{{\rm loc}}(z)$ have a natural structure of Lie algebras, and moreover $\mfrak{g} \otimes_\C \eus{P}^{\mfrak{g},\mfrak{p}}_{{\rm loc}}(z)$ and $\eus{C}^\mfrak{g}_{{\rm loc}}(z) \otimes_\C \eus{P}^{\mfrak{g},\mfrak{p}}_{{\rm loc}}(z)$ are $\mfrak{g} \otimes_\C \eus{P}^{\mfrak{g},\mfrak{p}}_{{\rm loc}}(z)$-modules through the adjoint action.

The Lie bracket on $\mfrak{g} \otimes_\C \eus{P}^{\mfrak{g},\mfrak{p}}_{{\rm loc}}(z)$ and $\eus{C}^\mfrak{g}_{{\rm loc}}(z) \otimes_\C \eus{P}^{\mfrak{g},\mfrak{p}}_{{\rm loc}}(z)$ is denoted by $[\cdot\,,\cdot]$. Let us note that there is no possibility for confusion with the Lie bracket given by \eqref{eq:comm relations series}. On the one hand we have $[a(z),b(z)]=[a,b](z)$, on the other hand we get $[a(z),b(w)]=[a,b](w) \delta(z-w) + (a,b)c\,\partial_w \delta(z-w)$ for $a,b \in \mfrak{g}$, other interpretation does not make sense.
\medskip

\lemma{\label{lem:ad action power}
We have
\begin{align}
(\ad(u(x)))^k(a(z))=(\ad(u(z)))^k(a(z))+(-1)^{k-1}((\ad(u(z)))^{k-1}(\partial_z u(z)),a)c \label{eq:ad action power}
\end{align}
for all $a \in \mfrak{g}$ and $k\in \N$.}

\proof{We prove it by induction on $k \in \N$. For $k=1$ we may write
\begin{align*}
  \ad(u(x))(a(z))&=\sum_{m\in\Z}  \sum_{\alpha \in \Delta(\mfrak{u})} \sum_{n\in \Z} x_{\alpha,n}[f_{\alpha,n},a_m] z^{-m-1} \\
  &= \sum_{m\in\Z}  \sum_{\alpha \in \Delta(\mfrak{u})} \sum_{n\in \Z} x_{\alpha,n}[f_\alpha,a]_{n+m} z^{-m-1} +  \sum_{m\in\Z}  \sum_{\alpha \in \Delta(\mfrak{u})} \sum_{n\in \Z} n\delta_{n,-m} x_{\alpha,n}(f_\alpha,a)c z^{-m-1} \\
  &= \sum_{m\in\Z}  \sum_{\alpha \in \Delta(\mfrak{u})} \sum_{n\in \Z} x_{\alpha,n}z^n[f_\alpha,a]_{n+m} z^{-n-m-1} +   \sum_{\alpha \in \Delta(\mfrak{u})} \sum_{n\in \Z} n x_{\alpha,n}(f_\alpha,a)c z^{n-1} \\
  &=   \sum_{\alpha \in \Delta(\mfrak{u})} a^*_\alpha(z)[f_\alpha,a](z) +   \sum_{\alpha \in \Delta(\mfrak{u})}  \partial_za^*(z)(f_\alpha,a)c\\
  &=  \ad(u(z))(a(z)) + (\partial_z u(z),a)c
\end{align*}
for all $a \in \mfrak{g}$. Now, let us assume that it holds for $k$. Since $(\ad(u(z)))^k(a(z)) \in \eus{C}^\mfrak{g}_{{\rm loc}}(z) \otimes_\C \eus{P}^{\mfrak{g},\mfrak{p}}_{{\rm loc}}(z)$, we get
\begin{align*}
  (\ad(u(z)))^k(a(z))=\sum_{r=1}^{\dim \mfrak{g}} P_r(z)b_r(z),
\end{align*}
where $\{b_r;\, r=1,2,\dots,\dim \mfrak{g}\}$ is a basis of $\mfrak{g}$ and $P_r(z) \in \eus{P}^{\mfrak{g},\mfrak{p}}_{{\rm loc}}(z)$ for $r=1,2,\dots,\dim \mfrak{g}$. Then we have
\begin{align*}
  (\ad(u(x)))^{k+1}(a(z))&=[u(x),(\ad(u(x)))^k(a(z))]=[u(x),(\ad(u(z)))^k(a(z))] \\
  &=\sum_{r=1}^{\dim\mfrak{g}} \sum_{\alpha \in \Delta(\mfrak{u})} \sum_{n \in \Z} x_{\alpha,n}P_r(z)[f_{\alpha,n},b_r(z)] \\
  & = \sum_{m\in \Z}\sum_{r=1}^{\dim\mfrak{g}} \sum_{\alpha \in \Delta(\mfrak{u})} \sum_{n \in \Z} x_{\alpha,n}P_r(z)[f_{\alpha,n},b_{r,m}]z^{-m-1} \\
  &=\sum_{m \in \Z}\sum_{r=1}^{\dim\mfrak{g}} \sum_{\alpha \in \Delta(\mfrak{u})} \sum_{n \in \Z} x_{\alpha,n}P_r(z)\big([f_\alpha,b_r]_{n+m}+n \delta_{n,-m}(f_\alpha,b_r)c \big)z^{-m-1}\\
  &=\sum_{m \in \Z}\sum_{r=1}^{\dim\mfrak{g}} \sum_{\alpha \in \Delta(\mfrak{u})} \sum_{n \in \Z} x_{\alpha,n}z^nP_r(z)[f_\alpha,b_r]_{n+m}z^{-n-m-1}\\
   &\quad + \sum_{r=1}^{\dim\mfrak{g}} \sum_{\alpha \in \Delta(\mfrak{u})} \sum_{n \in \Z} n x_{\alpha,n}z^{n-1}P_r(z)(f_\alpha,b_r)c \\
  &=\sum_{r=1}^{\dim\mfrak{g}} \sum_{\alpha \in \Delta(\mfrak{u})} a^*_\alpha(z)P_r(z)[f_\alpha,b_r](z) + \sum_{r=1}^{\dim\mfrak{g}} \sum_{\alpha \in \Delta(\mfrak{u})} \partial_za^*_\alpha(z)P_r(z)(f_\alpha,b_r)c.
\end{align*}
Therefore, we obtain
\begin{align*}
  (\ad(u(x)))^{k+1}(a(z))&=(\ad(u(z)))^{k+1}(a(z))+(\partial_z u(z),(\ad(u(z)))^k(a))c \\
  &=(\ad(u(z)))^{k+1}(a(z))+(-1)^k((\ad(u(z)))^k(\partial_z u(z)),a)c,
\end{align*}
where use used that $(\cdot\,,\cdot)$ is a $\mfrak{g}$-invariant symmetric bilinear form on $\mfrak{g}$. Thus, we are done.}

\lemma{\label{lem:exp ad action}
We have
\begin{align}
  e^{-\ad(u(x))}a(z)= e^{-\ad(u(z))}a(z)- \bigg({e^{\ad(u(z))}-\id \over \ad(u(z))}\,\partial_zu(z),a\!\bigg)c
\end{align}
for all $a \in \mfrak{g}$.}

\proof{By \eqref{eq:ad action power}, we may write
\begin{align*}
  e^{-\ad(u(x))}a(z)&= \sum_{k=0}^\infty {(-1)^k (\ad(u(x)))^k \over k!}\,a(z) \\
  &=\sum_{k=0}^\infty {(-1)^k (\ad(u(z)))^k \over k!}\,a(z) - \sum_{k=1}^\infty {((\ad(u(z)))^{k-1}(\partial_zu(z)),a) \over k!}\,c \\
  &=e^{-\ad(u(z))} a(z) - \bigg({e^{\ad(u(z))}-\id \over \ad(u(z))}\,\partial_zu(z),a\!\bigg)c.
\end{align*}
The proof is complete.}

Now, we may rewrite Theorem \ref{thm:realization coordinates} into an equivalent and more compact form. From \eqref{eq:explicit realization coordinates} we have
\begin{align*}
  \sum_{m \in \Z} \pi(a_m)z^{-m-1} &= - \sum_{m \in \Z} \sum_{\alpha \in \Delta(\mfrak{u})} \sum_{n \in \Z} \partial_{x_{\alpha,n}} \bigg[{\ad(u(x))e^{\ad(u(x))} \over e^{\ad(u(x))}-{\rm id}_{\widebar{\mfrak{u}}_{{\rm nat}}}}\,(e^{-\ad(u(x))}a_m)_{\widebar{\mfrak{u}}_{{\rm nat}}}\bigg]_{\alpha,n} z^{-m-1}\\
   &\quad + \sum_{m \in \Z} (e^{-\ad(u(x))}a_m)_{\mfrak{p}_{\rm nat}}  z^{-m-1} \\
  &= - \sum_{\alpha \in \Delta(\mfrak{u})} \sum_{n \in \Z} \partial_{x_{\alpha,n}} \bigg[{\ad(u(x))e^{\ad(u(x))} \over e^{\ad(u(x))}-{\rm id}_{\widebar{\mfrak{u}}_{{\rm nat}}}}\,(e^{-\ad(u(x))}a(z))_{\widebar{\mfrak{u}}_{{\rm nat}}}\bigg]_{\alpha,n}\\
   &\quad + (e^{-\ad(u(x))}a(z))_{\mfrak{p}_{\rm nat}}.
\end{align*}
By Lemma \ref{lem:exp ad action}, we get
\begin{align}
   (e^{-\ad(u(x))}a(z))_{\mfrak{p}_{\rm nat}} = (e^{-\ad(u(z))} a(z))_\mfrak{p} - \bigg({e^{\ad(u(z))}-\id \over \ad(u(z))}\,\partial_zu(z),a\!\bigg)c
\end{align}
and
\begin{align}
  (e^{-\ad(u(x))}a(z))_{\widebar{\mfrak{u}}_{\rm nat}}=(e^{-\ad(u(z))}a(z))_{\widebar{\mfrak{u}}}.
\end{align}
By a similar computation as in the proof of Lemma \ref{lem:ad action power}, we obtain
\begin{align*}
  {\ad(u(x))e^{\ad(u(x))} \over e^{\ad(u(x))}-{\rm id}_{\widebar{\mfrak{u}}_{{\rm nat}}}}\,(e^{-\ad(u(x))}a(z))_{\widebar{\mfrak{u}}_{{\rm nat}}} = {\ad(u(z))e^{\ad(u(z))} \over e^{\ad(u(z))}-\id}\,(e^{-\ad(u(z))}a(z))_{\widebar{\mfrak{u}}},
\end{align*}
since $(e^{-\ad(u(z))}a(z))_{\widebar{\mfrak{u}}} \in \eus{C}^{\widebar{\mfrak{u}}}_{{\rm loc}}(z) \otimes_\C \eus{P}^{\mfrak{g},\mfrak{p}}_{{\rm loc}}(z)$. Because 
\begin{align}
{\ad(u(z))e^{\ad(u(z))} \over e^{\ad(u(z))}-\id}\,(e^{-\ad(u(z))}a(z))_{\widebar{\mfrak{u}}} \in \eus{C}^{\widebar{\mfrak{u}}}_{{\rm loc}}(z) \otimes_\C \eus{P}^{\mfrak{g},\mfrak{p}}_{{\rm loc}}(z),
\end{align}
we have
\begin{align}
  {\ad(u(z))e^{\ad(u(z))} \over e^{\ad(u(z))}-\id}\,(e^{-\ad(u(z))}a(z))_{\widebar{\mfrak{u}}} = \sum_{\alpha \in \Delta(\mfrak{u})} P_\alpha(z) f_\alpha(z),
\end{align}
where $P_\alpha(z) \in \eus{P}^{\mfrak{g},\mfrak{p}}_{{\rm loc}}(z)$ and $\{f_\alpha;\, \alpha \in \Delta(\mfrak{u})\}$ is a basis of $\widebar{\mfrak{u}}$. Thus, we may write
\begin{multline*}
 \sum_{\alpha \in \Delta(\mfrak{u})} \sum_{n \in \Z} \partial_{x_{\alpha,n}} \bigg[{\ad(u(x))e^{\ad(u(x))} \over e^{\ad(u(x))}-{\rm id}_{\widebar{\mfrak{u}}_{{\rm nat}}}}\,(e^{-\ad(u(x))}a(z))_{\widebar{\mfrak{u}}_{{\rm nat}}}\bigg]_{\alpha,n}\\ = \sum_{\alpha \in \Delta(\mfrak{u})} \sum_{n \in \Z} \sum_{\beta \in \Delta(\mfrak{u})} \partial_{x_{\alpha,n}} P_\beta(z) [f_\beta(z)]_{\alpha,n} = \sum_{\alpha \in \Delta(\mfrak{u})} \sum_{n \in \Z} \partial_{x_{\alpha,n}} z^{-n-1} P_\alpha(z)\\ = \sum_{\alpha \in \Delta(\mfrak{u})} a_\alpha(z)P_\alpha(z) = \sum_{\alpha\in \Delta(\mfrak{u})} a_\alpha(z) \bigg[{\ad(u(z)) e^{\ad(u(z))} \over e^{\ad(u(z))}- \id_{\widebar{\mfrak{u}}}}\, (e^{-\ad(u(z))}a)_{\widebar{\mfrak{u}}}\bigg]_\alpha,
\end{multline*}
where we used 
\begin{align}
  {\ad(u(z))e^{\ad(u(z))} \over e^{\ad(u(z))}-\id}\,(e^{-\ad(u(z))}a)_{\widebar{\mfrak{u}}} = \sum_{\alpha \in \Delta(\mfrak{u})} P_\alpha(z) f_\alpha.
\end{align}
Taken altogether, we easily obtain the following theorem obviously equivalent to Theorem \ref{thm:realization coordinates}.
\medskip

\theorem{\label{thm:realization power series}
The embedding of $\widehat{\mfrak{g}}$ into $\eus{R}^{\mfrak{g},\mfrak{p}}_{{\rm loc}}$ is given by $\pi(c)=c$ and
\begin{multline}
  \pi(a(z))= - \sum_{\alpha\in \Delta(\mfrak{u})} a_\alpha(z) \bigg[{\ad(u(z)) e^{\ad(u(z))} \over e^{\ad(u(z))}- \id}\, (e^{-\ad(u(z))}a)_{\widebar{\mfrak{u}}}\bigg]_\alpha \\ + (e^{-\ad(u(z))} a(z))_\mfrak{p} - \bigg({e^{\ad(u(z))}-\id \over \ad(u(z))}\,\partial_z u(z),a\!\bigg)c,  \label{eq:explicit realization series}
\end{multline}
for all $a\in \mfrak{g}$, where
\begin{align}
   u(z) = \sum_{\alpha \in \Delta(\mfrak{u})} a_\alpha^*(z)f_\alpha
\end{align}
and
$[b]_\alpha$ is the $\alpha$-th coordinate of $b \in \widebar{\mfrak{u}}$ with respect to the basis $\{f_\alpha;\, \alpha \in \Delta(\mfrak{u})\}$ of $\widebar{\mfrak{u}}$. In particular, we have
\begin{align}
  \pi(a(z))=-\sum_{\alpha\in \Delta(\mfrak{u})} a_\alpha(z)\bigg[{\ad(u(z)) \over e^{\ad(u(z))}-\id}\,a\bigg]_\alpha  \label{eq:action power series op nilradical}
\end{align}
for $a \in \widebar{\mfrak{u}}$ and
\begin{align}
   \pi(a(z))=  \sum_{\alpha\in \Delta(\mfrak{u})} a_\alpha(z) [\ad(u(z))(a)]_\alpha  + a(z) \label{eq:action power series levi}
\end{align}
for $a \in \mfrak{l}$.}

\proof{If we introduce for greater clarity
\begin{align}
  D(a,z)&= - \sum_{\alpha\in \Delta(\mfrak{u})} a_\alpha(z) \bigg[{\ad(u(z)) e^{\ad(u(z))} \over e^{\ad(u(z))}- \id}\, (e^{-\ad(u(z))}a)_{\widebar{\mfrak{u}}}\bigg]_\alpha, \\
  A(a,z)&= (e^{-\ad(u(z))} a(z))_\mfrak{p}, \\
  C(a,z)&= - \bigg({e^{\ad(u(z))}-\id \over \ad(u(z))}\, \partial_z u(z),a\!\bigg)c
\end{align}
for $a \in \mfrak{g}$, then
\begin{align}
  \pi(a(z))=D(a,z)+A(a,z)+C(a,z).
\end{align}
By \eqref{eq:comm relations series}, we have
\begin{align*}
  \pi([a(z),b(w)])&=\pi([a,b](w))\delta(z-w)+(a,b)c\,\partial_w \delta(z-w) \\
  &=\big(D([a,b],w)+A([a,b],w)+C([a,b],w)\big)\delta(z-w)+ (a,b)c\,\partial_w \delta(z-w)
\end{align*}
for all $a,b \in \mfrak{g}$. On the other hand, we get
\begin{align*}
  [\pi(a(z)),\pi(b(w))]&=[D(a,z)+A(a,z)+C(a,z),D(b,w)+A(b,w)+C(b,w)] \\
  &= [D(a,z),D(b,w)]+[D(a,z),A(b,w)]+[A(a,z),D(b,w)] \\
  & \quad + [D(a,z),C(b,w)] + [C(a,z),D(b,w)] + [A(a,z),A(b,w)],
\end{align*}
where we used that $[A(a,z),C(b,w)]=0$ and $[C(a,z),A(b,w)]=0$. In addition, we have
\begin{align*}
  (e^{-\ad(u(z))}a(z))_\mfrak{p}= \sum_{r=1}^{\dim \mfrak{p}} P_r^a(z)d_r(z) \qquad \text{and} \qquad  (e^{-\ad(u(z))}b(z))_\mfrak{p}= \sum_{r=1}^{\dim \mfrak{p}} P_r^b(z)d_r(z)
\end{align*}
where $\{d_r;\, r=1,2,\dots, \dim \mfrak{p}\}$ is a basis of $\mfrak{p}$ and
$P_r^a(z), P_r^b(z) \in \eus{P}^{\mfrak{g},\mfrak{p}}_{{\rm loc}}(z)$ for $r=1,2,\dots, \dim \mfrak{p}$, hence we may write
\begin{align*}
  [A(a,z),A(b,w)]& =[(e^{-\ad(u(z))}a(z))_{\mfrak{p}},(e^{-\ad(u(w))}b(w))_{\mfrak{p}}] = \sum_{r,s=1}^{\dim \mfrak{p}} P_r^a(z) P_s^b(w)[d_r(z),d_s(w)] \\
  & = \sum_{r,s=1}^{\dim \mfrak{p}} P_r^a(z) P_s^b(w) \big([d_r,d_s](w)\delta(z-w) + (d_r,d_s)c\, \partial_w \delta(z-w)\big) \\
  & = \sum_{r,s=1}^{\dim \mfrak{p}} P_r^a(w) P_s^b(w) [d_r,d_s](w)\delta(z-w) + \sum_{r,s=1}^{\dim \mfrak{p}} P_r^a(z) P_s^b(w) (d_r,d_s)c\, \partial_w \delta(z-w) \\
  & = [(e^{-\ad(u(z))}a(w))_{\mfrak{p}},(e^{-\ad(u(w))}b(w))_{\mfrak{p}}] \delta(z-w) \\
   & \quad + ((e^{-\ad(u(z))}a)_{\mfrak{p}},(e^{-\ad(u(w))}b)_{\mfrak{p}}) c\,\partial_w \delta(z-w),
\end{align*}
where we used \eqref{eq:comm relations series} in the third equality.
Thus, proving that $\pi \colon \widehat{\mfrak{g}} \rarr \eus{R}^{\mfrak{g},\mfrak{p}}_{{\rm loc}}$ is a homomorphism of Lie algebras is equivalent to the following system of equations
\begin{gather}
  [D(a,z),D(b,w)]=D([a,b],w)\delta(z-w), \\
  [D(a,z),A(b,w)]+[A(a,z),D(b,w)]+ [A(a,z),A(b,w)]_{{\rm n}}=A([a,b],w)\delta(z-w),\\
  \begin{split}
  [D(a,z),C(b,w)] + [C(a,z),D(b,w)]+[A(a,z),A(b,w)]_{{\rm c}}& \\
   &\hspace{-0.13\columnwidth}=C([a,b],w)\delta(z-w)+ (a,b)c\,\partial_w \delta(z-w)
  \end{split}
\end{gather}
for all $a, b \in \mfrak{g}$, where
\begin{align}
[A(a,z),A(b,w)]_{{\rm c}} &= ((e^{-\ad(u(z))}a)_{\mfrak{p}},(e^{-\ad(u(w))}b)_{\mfrak{p}}) c\,\partial_w \delta(z-w), \label{eq:comm A and A central} \\
[A(a,z),A(b,w)]_{{\rm n}} &= [(e^{-\ad(u(z))}a(w))_{\mfrak{p}},(e^{-\ad(u(w))}b(w))_{\mfrak{p}}] \delta(z-w). \label{eq:comm A and A noncentral}
\end{align}
A proof of the previous system of equations is a subject of the following lemmas, which then completes the proof of the present theorem.}

Let us denote
\begin{align}
 g(\ad(u(z)))= {e^{\ad(u(z))}-\id \over \ad(u(z))}
\end{align}
and
\begin{align}
  T(a,z)={\ad(u(z)) e^{\ad(u(z))} \over e^{\ad(u(z))}- \id}\, (e^{-\ad(u(z))}a)_{\widebar{\mfrak{u}}} \label{eq:T(a,z) def}
\end{align}
for $a \in \mfrak{g}$. Then we have
\begin{align}
  T(a,z) = \sum_{\alpha \in \Delta(\mfrak{u})}  T_\alpha(a,z) f_\alpha,
\end{align}
where $T_\alpha(a,z)=[T(a,z)]_\alpha$ for $\alpha \in \Delta(\mfrak{u})$.
\medskip

\lemma{\label{lem:commutator a(z) and u(w)}
We have
\begin{align}
  [a_\alpha(z),u(w)]=f_\alpha \delta(z-w)\qquad \text{and} \qquad [a_\alpha(z),\partial_wu(w)] = f_\alpha \partial_w \delta(z-w)  \label{eq:commutator a(z) and u(w)}
\end{align}
for all $\alpha \in \Delta(\mfrak{u})$. Moreover, we have
\begin{align}
\begin{gathered}[]
  [D(a,z),u(w)]=-T(a,w)\delta(z-w), \\
  [D(a,z),\partial_wu(w)]=-T(a,w)\partial_w \delta(z-w) - \partial_w T(a,w) \delta(z-w) \label{eq:commutator D(a,z) and u(w)}
\end{gathered}
\end{align}
for all $a \in \mfrak{g}$.}

\proof{By \eqref{eq:u(z) def}, we have
\begin{align*}
  [a_\alpha(z),u(w)]= \sum_{\beta \in \Delta(\mfrak{u})} [a_\alpha(z),a^*_\beta(w)]f_\beta=  \sum_{\beta \in \Delta(\mfrak{u})} \delta_{\alpha\beta} f_\beta \delta(z-w)= f_\alpha \delta(z-w),
\end{align*}
therefore
\begin{align*}
  [a_\alpha(z),\partial_w u(w)]=\partial_w[a_\alpha(z),u(w)]=f_\alpha \partial_w \delta(z-w).
\end{align*}
Furthermore, we may write
\begin{align*}
   [D(a,z),u(w)]&=- \sum_{\alpha \in \Delta(\mfrak{u})} T_\alpha(a,z)[a_\alpha(z),u(w)] = - \sum_{\alpha \in \Delta(\mfrak{u})} T_\alpha(a,z) f_\alpha \delta(z-w) \\
   &= -T(a,z)\delta(z-w) = -T(a,w)\delta(z-w)
\end{align*}
and similarly
\begin{align*}
   [D(a,z),\partial_w u(w)]&=- \sum_{\alpha \in \Delta(\mfrak{u})} T_\alpha(a,z)[a_\alpha(z),\partial_w u(w)] = - \sum_{\alpha \in \Delta(\mfrak{u})} T_\alpha(a,z) f_\alpha \partial_w \delta(z-w) \\
   &= -T(a,z)\partial_w \delta(z-w) = -T(a,w)\partial_w \delta(z-w) - \partial_w T(a,w)\delta(z-w).
\end{align*}
We are done.}

\proposition{\label{prop:relation}
We have the following identities
\begin{enumerate}
\item[1)]
\begin{align}
  \bigg({{\rm d} \over {\rm d}t}_{|t=0} e^{\ad(u(z)+tx(z))}\!\bigg)e^{-\ad(u(z))}=\ad\!\bigg({e^{\ad(u(z))}-\id \over \ad(u(z))}\,x(z)\!\!\bigg) \label{eq:relation I}
\end{align}
for all $x(z) \in \eus{C}^\mfrak{g}_{{\rm loc}}(z) \otimes_\C \eus{P}^{\mfrak{g},\mfrak{p}}_{{\rm loc}}(z)$,
\item[2)]
\begin{multline}
   \bigg[{e^{\ad(u(z))}- \id \over \ad(u(z))}\,x(z),{e^{\ad(u(z))}- \id \over \ad(u(z))}\,y(z)\bigg] \\
   = {{\rm d} \over {\rm d}t}_{|t=0} {e^{\ad(u(z)+tx(z))}-\id \over \ad(u(z)+tx(z))}\,y(z) - {{\rm d} \over {\rm d}t}_{|t=0} {e^{\ad(u(z)+ty(z))}-\id \over \ad(u(z)+ty(z))}\,x(z) \label{eq:relation II}
\end{multline}
for all $x(z), y(z) \in \eus{C}^\mfrak{g}_{{\rm loc}}(z) \otimes_\C \eus{P}^{\mfrak{g},\mfrak{p}}_{{\rm loc}}(z)$.
\end{enumerate}}

\proof{First of all, for $x(z), y(z) \in \eus{C}^\mfrak{g}_{{\rm loc}}(z) \otimes_\C \eus{P}^{\mfrak{g},\mfrak{p}}_{{\rm loc}}(z)$ we have
\begin{align*}
  {{\rm d} \over {\rm d}t}_{|t=0} (\ad(u(z)+tx(z)))^n(y(z))&= \sum_{k=0}^{n-1}\, (\ad(u(z)))^k(\ad(x(z)))(\ad(u(z)))^{n-1-k}(y(z))\\
  &= \sum_{k=0}^{n-1}\, (\ad(u(z)))^k([x(z),(\ad(u(z)))^{n-1-k}(y(z))]) \\
  & = \sum_{k=0}^{n-1} \sum_{j=0}^k \binom{k}{j} [(\ad(u(z)))^j(x(z)),(\ad(u(z)))^{n-1-j}(y(z))] \\
  & = \sum_{j=0}^{n-1} \sum_{k=j}^{n-1} \binom{k}{j} [(\ad(u(z)))^j(x(z)),(\ad(u(z)))^{n-1-j}(y(z))] \\
  & = \sum_{j=0}^{n-1} \!\binom{n}{j+1} [(\ad(u(z)))^j(x(z)),(\ad(u(z)))^{n-1-j}(y(z))],
\end{align*}
where we used the Leibniz rule for the derivation $\ad(u(z))$ of the Lie algebra $\eus{C}^\mfrak{g}_{{\rm loc}}(z) \otimes_\C \eus{P}^{\mfrak{g},\mfrak{p}}_{{\rm loc}}(z)$ in the third equality and the binomial identity $\sum_{k=j}^{n-1} \!\binom{k}{j} =\binom{n}{j+1}$ in the last equality. Therefore
\begin{align}
  {{\rm d} \over {\rm d}t}_{|t=0} (\ad(u(z)+tx(z)))^n(y(z))= \sum_{j=0}^{n-1} \!\binom{n}{j+1} [(\ad(u(z)))^j(x(z)),(\ad(u(z)))^{n-1-j}(y(z))] \label{eq:derivative ad}
\end{align}
for all $x(z), y(z) \in \eus{C}^\mfrak{g}_{{\rm loc}}(z) \otimes_\C \eus{P}^{\mfrak{g},\mfrak{p}}_{{\rm loc}}(z)$.

Now, we may write
\begin{align*}
\begin{split}
  {{\rm d} \over {\rm d}t}_{|t=0} e^{\ad (u(z)+tx(z))}y(z)&= {{\rm d} \over {\rm d}t}_{|t=0} \sum_{n=0}^\infty {(\ad (u(z)+tx(z)))^n \over n!}\, y(z) \\
   & \hspace{-0.02\columnwidth} = \sum_{n=1}^\infty {1 \over n!} \sum_{j=0}^{n-1}\! \binom{n}{j+1} [(\ad(u(z)))^j(x(z)), (\ad(u(z)))^{n-1-j}(y(z))]\\
   & \hspace{-0.02\columnwidth} =\sum_{j=0}^\infty  \sum_{n=j+1}^\infty  {1 \over (j+1)!} {1 \over (n-1-j)!}\,  [(\ad(u(z)))^j(x(z)), (\ad(u(z)))^{n-1-j}(y(z))] \\
   & \hspace{-0.02\columnwidth} = \bigg[{e^{\ad(u(z))} - \id \over \ad(u(z))}\,x(z), e^{\ad(u(z))}y(z)\bigg] \\
   & \hspace{-0.02\columnwidth} = \ad\!\bigg({e^{\ad(u(z))} - \id \over \ad(u(z))}\,x(z)\!\!\bigg)(e^{\ad(u(z))}y(z)),
\end{split}
\end{align*}
where we used \eqref{eq:derivative ad}. If we take $e^{-\ad(u(z))}y(z)$ instead of $y(z)$, we obtain the required statement.

For the second identity, we may write
\begin{align*}
\begin{split}
  {{\rm d} \over {\rm d}t}_{|t=0} {e^{\ad(u(z)+tx(z))} - \id \over \ad(u(z)+tx(z))}\, y(z) & = {{\rm d} \over {\rm d}t}_{|t=0} \sum_{n=0}^\infty {(\ad(u(z)+tx(z)))^n \over (n+1)!}\,y(z) \\
  & \hspace{-0.02\columnwidth} = \sum_{n=1}^\infty {1 \over (n+1)!} \sum_{j=0}^{n-1}\! \binom{n}{j+1} [(\ad(u(z)))^j(x(z)), (\ad(u(z)))^{n-1-j}(y(z))]
\end{split}
\end{align*}
and
\begin{align*}
\begin{split}
  {{\rm d} \over {\rm d}t}_{|t=0} {e^{\ad(u(z)+ty(z))} - \id \over \ad(u(z)+ty(z))}\, x(z) & =  {{\rm d} \over {\rm d}t}_{|t=0} \sum_{n=0}^\infty {(\ad(u(z)+ty(z)))^n \over (n+1)!}\, x(z) \\
  & \hspace{-0.02\columnwidth} = \sum_{n=1}^\infty {1 \over (n+1)!} \sum_{j=0}^{n-1} \!\binom{n}{j+1} [(\ad(u(z)))^j(y(z)), (\ad(u(z)))^{n-1-j}(x(z))] \\
  & \hspace{-0.02\columnwidth} =- \sum_{n=1}^\infty {1 \over (n+1)!} \sum_{j=0}^{n-1} \binom{n}{j} [(\ad(u(z)))^j(x(z)), (\ad(u(z)))^{n-1-j}(y(z))].
\end{split}
\end{align*}
Hence, for the right hand side of \eqref{eq:relation II} we have
\begin{align*}
\begin{split}
  \sum_{n=1}^\infty {1 \over (n+1)!} \sum_{j=0}^{n-1} \bigg(\!\!\binom{n}{j} + \binom{n}{j+1}\!\!\bigg) [(\ad(u(z)))^j(x(z)), (\ad(u(z)))^{n-1-j}(y(z))]& \\
  &\hspace{-0.5\columnwidth} = \sum_{n=1}^\infty {1 \over (n+1)!} \sum_{j=0}^{n-1} \binom{n+1}{j+1} [(\ad(u(z)))^j(x(z)), (\ad(u(z)))^{n-1-j}(y(z))] \\
  &\hspace{-0.5\columnwidth} = \sum_{j=0}^\infty \sum_{n=j+1}^\infty {1 \over (j+1)!} {1 \over (n-j)!}  [(\ad(u(z)))^j(x(z)), (\ad(u(z)))^{n-1-j}(y(z))  \\
  &\hspace{-0.5\columnwidth} = \bigg[{e^{\ad(u(z))} - \id \over \ad(u(z))}\, x(z),  {e^{\ad(u(z))} - \id \over \ad(u(z))}\, y(z)\bigg],
\end{split}
\end{align*}
which completes the claim.}

\lemma{We have
\begin{multline}
  [D(a,z),C(b,w)] + [C(a,z),D(b,w)] + [A(a,z),A(b,w)]_{{\rm c}} \\
  =C([a,b],w)\delta(z-w)+ (a,b)c\,\partial_w \delta(z-w) \label{eq:comm D and C}
\end{multline}
for all $a,b \in \mfrak{g}$.}

\proof{We have
\begin{align*}
  [D(a,z),C(b,w)]
  & = -[D(a,z),(g(\ad(u(w)))(\partial_w u(w)),b)]c \\
  &=-([D(a,z),g(\ad(u(w)))](\partial_w u(w)),b)c  - (g(\ad(u(w)))([D(a,z),\partial_w u(w)]),b)c.
\end{align*}
The first summand may be written as
\begin{align*}
    ([D(a,z),g(\ad(u(w)))](\partial_w u(w)),b)c = -{{\rm d}\over {\rm d}t}_{|{t=0}} (g(\ad(u(w)+tT(a,w)))(\partial_w u(w)),b)c\,\delta(z-w),
\end{align*}
and the second simplifies to
\begin{align*}
\begin{split}
(g(\ad(u(w)))([D(a,z),\partial_w u(w)]),b)c &\\
& \hspace{-0.24\columnwidth} = -(g(\ad(u(w)))(T(a,w)),b)c\,\partial_w \delta(z-w) - (g(\ad(u(w)))(\partial_wT(a,w)),b)c\,\delta(z-w)  \\
& \hspace{-0.24\columnwidth} = -(g(\ad(u(w)))(T(a,w)),b)c\,\partial_w \delta(z-w) + (\partial_wg(\ad(u(w)))(T(a,w)),b)c\,\delta(z-w)  \\
& \hspace{-0.24\columnwidth} \quad - \partial_w(g(\ad(u(w)))(T(a,w)),b)c\, \delta(z-w),
\end{split}
\end{align*}
where we used \eqref{eq:commutator D(a,z) and u(w)}. Similarly, we obtain
\begin{align*}
  [C(a,z),D(b,w)]
  & = [D(b,w),(g(\ad(u(z)))(\partial_z u(z)),a)]c \\
  &= ([D(b,w),g(\ad(u(z)))](\partial_z u(z)),a)c  + (g(\ad(u(z)))([D(b,w),\partial_z u(z)]),a)c
\end{align*}
with
\begin{align*}
    ([D(b,w),g(\ad(u(z)))](\partial_z u(z)),a)c = -{{\rm d}\over {\rm d}t}_{|{t=0}} (g(\ad(u(w)+tT(b,w)))(\partial_w u(w)),a)c\,\delta(z-w),
\end{align*}
and
\begin{align*}
\begin{split}
 (g(\ad(u(z)))([D(b,w),\partial_z u(z)]),a)c &\\
 & \hspace{-0.22\columnwidth} = (g(\ad(u(z)))(T(b,w)),a)c\,\partial_w \delta(z-w) \\
 & \hspace{-0.22\columnwidth} = (g(\ad(u(w)))(T(b,w)),a)c\,\partial_w \delta(z-w)
   + (\partial_w g(\ad(u(w)))(T(b,w)),a)c\,\delta(z-w),
\end{split}
\end{align*}
where we used \eqref{eq:commutator D(a,z) and u(w)} and Proposition \ref{prop:delta function}. Finally, we have
\begin{align*}
  [A(a,z),A(b,w)]_{{\rm c}}&= ((e^{-\ad(u(z))}a)_{\mfrak{p}},(e^{-\ad(u(w))}b)_{\mfrak{p}}) c\,\partial_w \delta(z-w) \\
  &=((e^{-\ad(u(w))}a)_{\mfrak{p}},(e^{-\ad(u(w))}b)_{\mfrak{p}}) c\,\partial_w \delta(z-w)\\
   &\quad + (\partial_w(e^{-\ad(u(w))}a)_{\mfrak{p}},(e^{-\ad(u(w))}b)_{\mfrak{p}}) c\,\delta(z-w).
\end{align*}
Collecting all terms together and comparing the coefficients in front of $c\,\partial_w \delta(z-w)$ and $c\,\delta(z-w)$ in \eqref{eq:comm D and C}, we obtain that \eqref{eq:comm D and C} is equivalent to
\begin{multline}
   ((e^{-\ad(u(w))}a)_{\widebar{\mfrak{u}}},(e^{-\ad(u(w))}b)_\mfrak{u}) + ((e^{-\ad(u(w))}a)_\mfrak{u},(e^{-\ad(u(w))}b)_{\widebar{\mfrak{u}}}) \\ + ((e^{-\ad(u(w))}a)_{\mfrak{p}},(e^{-\ad(u(w))}b)_{\mfrak{p}})= (a,b)  \label{eq:derivative delta coeff}
\end{multline}
and
\begin{multline}
 {{\rm d}\over {\rm d}t}_{|{t=0}} (g(\ad(u(w)+tT(a,w)))(\partial_w u(w)),b) -
 {{\rm d}\over {\rm d}t}_{|{t=0}} (g(\ad(u(w)+tT(b,w)))(\partial_wu(w)),a) \\
 + (\partial_w g(\ad(u(w)))(T(b,w)),a) - (\partial_w g(\ad(u(w)))(T(a,w)),b)\\
 + \partial_w((e^{-\ad(u(w))}a)_{\widebar{\mfrak{u}}},(e^{-\ad(u(w))}b)_\mfrak{u}) +(\partial_w(e^{-\ad(u(w))}a)_{\mfrak{p}},(e^{-\ad(u(w))}b)_{\mfrak{p}}) \\
  = - (g(\ad(u(w)))(\partial_w u(w)),[a,b]) \label{eq:delta coeff}
\end{multline}
for all $a,b \in \mfrak{g}$, where we used
\begin{align*}
  (g(\ad(u(w)))(T(a,w)),b)=(e^{\ad(u(w))}(e^{-\ad(u(w))}a)_{\widebar{\mfrak{u}}},b)= ((e^{-\ad(u(w))}a)_{\widebar{\mfrak{u}}},(e^{-\ad(u(w))}b)_\mfrak{u})
\end{align*}
and
\begin{align*}
  (a,g(\ad(u(w)))(T(b,w)))=(a,e^{\ad(u(w))}(e^{-\ad(u(w))}b)_{\widebar{\mfrak{u}}})= ((e^{-\ad(u(w))}a)_\mfrak{u},(e^{-\ad(u(w))}b)_{\widebar{\mfrak{u}}}).
\end{align*}
Since we have
\begin{align*}
  ((e^{-\ad(u(w))}a)_{\mfrak{p}},(e^{-\ad(u(w))}b)_{\mfrak{p}}) = ((e^{-\ad(u(w))}a)_{\mfrak{l}},(e^{-\ad(u(w))}b)_{\mfrak{l}}),
\end{align*}
we immediately may rewrite the left hand side of \eqref{eq:derivative delta coeff} into the form
\begin{multline*}
   ((e^{-\ad(u(w))}a)_{\widebar{\mfrak{u}}},(e^{-\ad(u(w))}b)_\mfrak{u}) + ((e^{-\ad(u(w))}a)_\mfrak{u},(e^{-\ad(u(w))}b)_{\widebar{\mfrak{u}}}) \\ + ((e^{-\ad(u(w))}a)_{\mfrak{p}},(e^{-\ad(u(w))}b)_{\mfrak{p}}) = (e^{-\ad(u(w))}a,e^{-\ad(u(w))}b) = (a,b)
\end{multline*}
for all $a,b \in \mfrak{g}$. Therefore, we have proved \eqref{eq:derivative delta coeff}.

Moreover, by \eqref{eq:relation II} we have
\begin{align*}
\begin{split}
   {{\rm d} \over {\rm d}t}_{|t=0} {e^{\ad(u(w)+tT(a,w))}-\id \over \ad(u(w)+tT(a,w))}\, \partial_w u(w) - {{\rm d} \over {\rm d}t}_{|t=0} {e^{\ad(u(w)+t\partial_w u(w))}-\id \over \ad(u(w)+t\partial_w u(w))}\,T(a,w)&\\
     & \hspace{-0.43\columnwidth} =\bigg[{e^{\ad(u(w))}- \id \over \ad(u(w))}\, T(a,w),{e^{\ad(u(w))}- \id \over \ad(u(w))}\, \partial_w u(w) \bigg] \\
     & \hspace{-0.43\columnwidth} =-\ad\!\bigg({e^{\ad(u(w))}- \id \over \ad(u(w))}\,\partial_w u(w)\!\!\bigg)\!\big(e^{\ad(u(w))}(e^{-\ad(u(w))}a)_{\widebar{\mfrak{u}}}\big) \\
     & \hspace{-0.43\columnwidth} = -\bigg({{\rm d} \over {\rm d}t}_{|t=0} e^{\ad(u(w)+t\partial_w u(w))}\!\bigg)e^{-\ad(u(w))}e^{\ad(u(w))}(e^{-\ad(u(w))}a)_{\widebar{\mfrak{u}}} \\
     & \hspace{-0.43\columnwidth} =-(\partial_w e^{\ad(u(w))})(e^{-\ad(u(w))}a)_{\widebar{\mfrak{u}}}
\end{split}
\end{align*}
for all $a \in \mfrak{g}$. This gives us
\begin{multline*}
  {{\rm d}\over {\rm d}t}_{|{t=0}} (g(\ad(u(w)+tT(a,w)))(\partial_w u(w)),b) - {{\rm d}\over {\rm d}t}_{|{t=0}} (g(\ad(u(w)+t\partial_w u(w)))(T(a,w)),b)  \\
  =  -((\partial_w e^{\ad(u(w))})(e^{-\ad(u(w))}a)_{\widebar{\mfrak{u}}},b) = -((e^{-\ad(u(w))}a)_{\widebar{\mfrak{u}}},\partial_w(e^{-\ad(u(w))}b)_\mfrak{u})
\end{multline*}
for all $a,b \in \mfrak{g}$, where we used
\begin{align*}
((\partial_w e^{\ad(u(w))})(e^{-\ad(u(w))}a)_{\widebar{\mfrak{u}}},b)= ((e^{-\ad(u(w))}a)_{\widebar{\mfrak{u}}},\partial_w(e^{-\ad(u(w))}b)_\mfrak{u}),
\end{align*}
which follows from the following computation
\begin{align*}
\begin{split}
  ((\partial_w e^{\ad(u(w))})(e^{-\ad(u(w))}a)_{\widebar{\mfrak{u}}},b)&\\
  & \hspace{-0.09\columnwidth} =\partial_w (e^{\ad(u(w))}(e^{-\ad(u(w))}a)_{\widebar{\mfrak{u}}},b) -
  (e^{\ad(u(w))}\partial_w(e^{-\ad(u(w))}a)_{\widebar{\mfrak{u}}},b) \\
  & \hspace{-0.09\columnwidth} =\partial_w ((e^{-\ad(u(w))}a)_{\widebar{\mfrak{u}}},(e^{-\ad(u(w))}b)_\mfrak{u})  - (\partial_w(e^{-\ad(u(w))}a)_{\widebar{\mfrak{u}}},(e^{-\ad(u(w))}b)_\mfrak{u}) \\
  & \hspace{-0.09\columnwidth} = ((e^{-\ad(u(w))}a)_{\widebar{\mfrak{u}}},\partial_w(e^{-\ad(u(w))}b)_\mfrak{u}).
\end{split}
\end{align*}
Therefore, we have
\begin{multline*}
  {{\rm d}\over {\rm d}t}_{|{t=0}} (g(\ad(u(w)+tT(a,w)))(\partial_w u(w)),b) - {{\rm d}\over {\rm d}t}_{|{t=0}} (g(\ad(u(w)+t\partial_w u(w)))(T(a,w)),b)  \\
   = -((e^{-\ad(u(w))}a)_{\widebar{\mfrak{u}}},\partial_w(e^{-\ad(u(w))}b)_\mfrak{u})
\end{multline*}
and also
\begin{multline*}
  {{\rm d}\over {\rm d}t}_{|{t=0}} (g(\ad(u(w)+tT(b,w)))(\partial_w u(w)),a) - {{\rm d}\over {\rm d}t}_{|{t=0}} (g(\ad(u(w)+t\partial_w u(w)))(T(b,w)),a)  \\
   = -(\partial_w(e^{-\ad(u(w))}a)_\mfrak{u},(e^{-\ad(u(w))}b)_{\widebar{\mfrak{u}}})
\end{multline*}
for all $a,b \in \mfrak{g}$. Taking the difference of the two equations above and applying the resulting formula on the left hand side of \eqref{eq:delta coeff}, we may rewrite the left hand side of \eqref{eq:delta coeff} into the form
\begin{multline*}
  -((e^{-\ad(u(w))}a)_{\widebar{\mfrak{u}}},\partial_w(e^{-\ad(u(w))}b)_\mfrak{u})+ (\partial_w (e^{-\ad(u(w))}a)_\mfrak{u},(e^{-\ad(u(w))}b)_{\widebar{\mfrak{u}}}) \\
  +  \partial_w((e^{-\ad(u(w))}a)_{\widebar{\mfrak{u}}},(e^{-\ad(u(w))}b)_\mfrak{u}) +(\partial_w(e^{-\ad(u(w))}a)_{\mfrak{p}},(e^{-\ad(u(w))}b)_{\mfrak{p}}) \\
  = (\partial_w (e^{-\ad(u(w))}a)_\mfrak{u},(e^{-\ad(u(w))}b)_{\widebar{\mfrak{u}}})
   + (\partial_w (e^{-\ad(u(w))}a)_{\widebar{\mfrak{u}}},(e^{-\ad(u(w))}b)_\mfrak{u}) \\
   + (\partial_w(e^{-\ad(u(w))}a)_{\mfrak{l}},(e^{-\ad(u(w))}b)_{\mfrak{l}}) \\
  = (\partial_w e^{-\ad(u(w))}a,e^{-\ad(u(w))}b).
\end{multline*}
Moreover, using \eqref{eq:relation I} we finally get
\begin{align*}
  (\partial_w e^{-\ad(u(w))}a,e^{-\ad(u(w))}b)&=\bigg(\!\!\ad\!\bigg({e^{-\ad(u(w))}-\id \over \ad(u(w))}\, \partial_w u(w)\!\!\bigg)(e^{-\ad(u(w))}a),e^{-\ad(u(w))}b\!\bigg) \\
  &= \bigg(\bigg[{e^{-\ad(u(w))}-\id \over \ad(u(w))}\, \partial_w u(w), e^{-\ad(u(w))}a\bigg],e^{-\ad(u(w))}b\!\bigg) \\
  &= -\bigg(\bigg[{e^{\ad(u(w))}-\id \over \ad(u(w))}\, \partial_w u(w),a\bigg],b\!\bigg) \\
  &= -\bigg({e^{\ad(u(w))}-\id \over \ad(u(w))}\, \partial_w u(w),[a,b]\!\bigg)
\end{align*}
for all $a,b \in \mfrak{g}$. So, we have proved \eqref{eq:delta coeff}. The proof is complete.}

\lemma{We have
\begin{align}
  [D(a,z),A(b,w)]+[A(a,z),D(b,w)]+ [A(a,z),A(b,w)]_{{\rm n}}=A([a,b],w)\delta(z-w) \label{eq:comm D and A}
\end{align}
for all $a,b \in \mfrak{g}$.}

\proof{We have
\begin{align*}
  [D(a,z),A(b,w)]&=[D(a,z),(e^{-\ad(u(w))} b(w))_\mfrak{p}]=-{{\rm d}\over {\rm d}t}_{|t=0} (e^{-\ad(u(w)+tT(a,w))}b(w))_\mfrak{p}\,\delta(z-w) \\
  &= -\bigg(\!\!\ad\!\bigg({e^{-\ad(u(w))}-\id \over \ad(u(w))}\,T(a,w)\!\!\bigg) (e^{-\ad(u(w))}b(w))\!\!\bigg)_{\!\mfrak{p}} \delta(z-w) \\
  &= [(e^{-\ad(u(w))}a)_{\widebar{\mfrak{u}}},e^{-\ad(u(w))}b(w)]_\mfrak{p}\, \delta(z-w) \\
  &= [(e^{-\ad(u(w))}a(w))_{\widebar{\mfrak{u}}},(e^{-\ad(u(w))}b(w))_\mfrak{u}]_\mfrak{p}\, \delta(z-w)
\end{align*}
for all $a,b \in \mfrak{g}$, were we used $[D(a,z),u(w)]=-T(a,w) \delta(z-w)$ in the second equality, Proposition \ref{prop:relation} in the third equality, \eqref{eq:T(a,z) def} in the fourth equality, and the discussion in the paragraph preceding Lemma \ref{lem:ad action power} in the last equality. By analogous computation we get
\begin{align*}
  [A(a,z),D(b,w)] = [(e^{-\ad(u(w))}a(w))_\mfrak{u},(e^{-\ad(u(w))}b(w))_{\widebar{\mfrak{u}}}]_\mfrak{p}\, \delta(z-w)
\end{align*}
for all $a,b \in \mfrak{g}$. Therefore, for the left hand side of \eqref{eq:comm D and A} we may write
\begin{multline*}
  [(e^{-\ad(u(w))}a(w))_{\widebar{\mfrak{u}}},(e^{-\ad(u(w))}b(w))_\mfrak{u}]_\mfrak{p}\, \delta(z-w) \\
  \shoveleft{\hspace{0.15\columnwidth} + [(e^{-\ad(u(w))}a(w))_\mfrak{u},(e^{-\ad(u(w))}b(w))_{\widebar{\mfrak{u}}}]_\mfrak{p}\, \delta(z-w)}\\
  \shoveright{+ [(e^{-\ad(u(w))}a(w))_\mfrak{p},(e^{-\ad(u(w))}b(w))_\mfrak{p}]\, \delta(z-w)\hspace{0.13\columnwidth}} \\
  \shoveright{= [(e^{-\ad(u(w))}a(w), e^{-\ad(u(w))}b(w)]_\mfrak{p}\, \delta(z-w)} \\
  =  (e^{-\ad(u(w))}[a,b](w))_\mfrak{p}\, \delta(z-w),
\end{multline*}
which gives the statement.}

\lemma{We have
\begin{align}
   [D(a,z),D(b,w)]=D([a,b],w)\delta(z-w)
\end{align}
for all $a,b \in \mfrak{g}$.}

\proof{Since we have
\begin{align*}
  [D(a,z),D(b,w)]&=\sum_{\alpha,\beta \in \Delta(\mfrak{u})} [a_\alpha(z)T_\alpha(a,z), a_\beta(w)T_\beta(b,w)] \\
  &= -\sum_{\alpha,\beta \in \Delta(\mfrak{u})}\! \big(a_\alpha(z)[a_\beta(w),T_\alpha(a,z)]T_\beta(b,w) - a_\beta(w)[a_\alpha(z),T_\beta(b,w)]T_\alpha(a,z)\big) \\
  &= -  \sum_{\alpha,\beta \in \Delta(\mfrak{u})} a_\alpha(w)\big([a_\beta(z),T_\alpha(a,w)]T_\beta(b,z)- [a_\beta(z),T_\alpha(b,w)]T_\beta(a,z) \big),
\end{align*}
it is enough to prove that
\begin{align*}
  T_\alpha([a,b],w)\delta(z-w)= \sum_{\beta \in \Delta(\mfrak{u})}\!\big([a_\beta(z),T_\alpha(a,w)]T_\beta(b,z)- [a_\beta(z),T_\alpha(b,w)]T_\beta(a,z) \big)
\end{align*}
for all $a,b \in \mfrak{g}$ and $\alpha \in \Delta(\mfrak{u})$, or equivalently
\begin{align*}
  T([a,b],w)\delta(z-w)&= \sum_{\beta \in \Delta(\mfrak{u})}\!\big([a_\beta(z),T(a,w)]T_\beta(b,z)- [a_\beta(z),T(b,w)]T_\beta(a,z)\big) \\
  &= - \big([D(b,z),T(a,w)]- [D(a,z),T(b,w)]\big)
\end{align*}
for all $a,b \in \mfrak{g}$. Further, we may write
\begin{align*}
  [D(a,z),T(b,w)]
  = -{{\rm d} \over {\rm d}t}_{|t=0} {\ad(u(w)+tT(a,w)) e^{\ad(u(w)+tT(a,w))} \over e^{\ad(u(w)+tT(a,w))} - \id}\,(e^{-\ad(u(w)+tT(a,w))}b)_{\widebar{\mfrak{u}}}\, \delta(z-w),
\end{align*}
since $[D(a,z),u(w)]=-T(a,w)\delta(z-w)$. Therefore, we have
\begin{align*}
    [D(a,z),T(b,w)]&=   -{{\rm d} \over {\rm d}t}_{|t=0} {\ad(u(w)+tT(a,w)) \over e^{\ad(u(w)+tT(a,w))} - \id}\, e^{\ad(u(w))}(e^{-\ad(u(w))}b)_{\widebar{\mfrak{u}}}\, \delta(z-w) \\
    & \quad  - {\ad(u(w)) \over e^{\ad(u(w))} - \id}\, {{\rm d} \over {\rm d}t}_{|t=0} e^{\ad(u(w)+tT(a,w))}
    (e^{-\ad(u(w))}b)_{\widebar{\mfrak{u}}}\, \delta(z-w) \\
    & \quad   - {\ad(u(w))  e^{\ad(u(w))} \over e^{\ad(u(w))} - \id}\,
    {{\rm d} \over {\rm d}t}_{|t=0}(e^{-\ad(u(w)+tT(a,w))}b)_{\widebar{\mfrak{u}}}\, \delta(z-w)
\end{align*}
for all $a,b \in \mfrak{g}$. From \eqref{eq:relation I} we obtain
\begin{align*}
  {{\rm d} \over {\rm d}t}_{|t=0} e^{\ad(u(w)+tT(a,w))}
    (e^{-\ad(u(w))}b)_{\widebar{\mfrak{u}}} &=\ad\!\bigg({e^{\ad(u(w))}-\id \over \ad(u(w))}\,T(a,w)\!\!\bigg)\! \big(e^{\ad(u(w))}(e^{-\ad(u(w))}b)_{\widebar{\mfrak{u}}}\big) \\
    &= [e^{\ad(u(w))}(e^{-\ad(u(w))}a)_{\widebar{\mfrak{u}}}, e^{\ad(u(w))}(e^{-\ad(u(w))}b)_{\widebar{\mfrak{u}}}] \\
    &= e^{\ad(u(w))}[(e^{-\ad(u(w))}a)_{\widebar{\mfrak{u}}},(e^{-\ad(u(w))}b)_{\widebar{\mfrak{u}}}],
\end{align*}
and similarly we get
\begin{align*}
   {{\rm d} \over {\rm d}t}_{|t=0} e^{-\ad(u(w)+tT(a,w))}b &= \ad\!\bigg({e^{-\ad(u(w))}-\id \over \ad(u(w))}\, T(a,w)\!\!\bigg) (e^{-\ad(u(w))}b) \\
   &= - [(e^{-\ad(u(w))}a)_{\widebar{\mfrak{u}}},e^{-\ad(u(w))}b].
\end{align*}
Finally, we have
\begin{align*}
\begin{split}
  {{\rm d} \over {\rm d}t}_{|t=0} {\ad(u(w)+tT(a,w)) \over e^{\ad(u(w)+tT(a,w))} - \id}\, e^{\ad(u(w))}(e^{-\ad(u(w))}b)_{\widebar{\mfrak{u}}} & \\
  & \hspace{-0.15\columnwidth} = {{\rm d} \over {\rm d}t}_{|t=0} {\ad(u(w)+tT(a,w)) \over e^{\ad(u(w)+tT(a,w))} - \id}\, {e^{\ad(u(w))}-\id \over \ad(u(w))}\, T(b,w) \\
  & \hspace{-0.15\columnwidth} = - {\ad(u(w)) \over e^{\ad(u(w))} - \id}\, {{\rm d} \over {\rm d}t}_{|t=0} {e^{\ad(u(w)+tT(a,w))}-\id \over \ad(u(w)+tT(a,w))}\, T(b,w).
\end{split}
\end{align*}
If we put all together, we obtain
\begin{align*}
  [D(a,z),T(b,w)]&= {\ad(u(w)) \over e^{\ad(u(w))} - \id}\, {{\rm d} \over {\rm d}t}_{|t=0} {e^{\ad(u(w)+tT(a,w))}-\id \over \ad(u(w)+tT(a,w))}\, T(b,w)\, \delta(z-w) \\
  &\quad + {\ad(u(w)) e^{\ad(u(w))} \over e^{\ad(u(w))} - \id}\, [(e^{-\ad(u(w))}a)_{\widebar{\mfrak{u}}},e^{-\ad(u(w))}b]_{\widebar{\mfrak{u}}}  \,\delta(z-w)\\
  &\quad - {\ad(u(w)) e^{\ad(u(w))} \over e^{\ad(u(w))} - \id}\, [(e^{-\ad(u(w))}a)_{\widebar{\mfrak{u}}},(e^{-\ad(u(w))}b)_{\widebar{\mfrak{u}}}] \,\delta(z-w),
\end{align*}
and therefore also
\begin{align*}
  [D(b,z),T(a,w)]&= {\ad(u(w)) \over e^{\ad(u(w))} - \id}\, {{\rm d} \over {\rm d}t}_{|t=0} {e^{\ad(u(w)+tT(b,w))}-\id \over \ad(u(w)+tT(b,w))}\, T(a,w)\, \delta(z-w) \\
  &\quad + {\ad(u(w)) e^{\ad(u(w))} \over e^{\ad(u(w))} - \id}\, [(e^{-\ad(u(w))}b)_{\widebar{\mfrak{u}}},e^{-\ad(u(w))}a]_{\widebar{\mfrak{u}}}  \,\delta(z-w)\\
  &\quad - {\ad(u(w)) e^{\ad(u(w))} \over e^{\ad(u(w))} - \id}\, [(e^{-\ad(u(w))}b)_{\widebar{\mfrak{u}}},(e^{-\ad(u(w))}a)_{\widebar{\mfrak{u}}}] \,\delta(z-w).
\end{align*}
Taking the difference of the two equations above, we get for the right hand side
\begin{multline*}
  {\ad(u(w)) \over e^{\ad(u(w))} - \id}\,\bigg({{\rm d} \over {\rm d}t}_{|t=0} {e^{\ad(u(w)+tT(a,w))}-\id \over \ad(u(w)+tT(a,w))}\, T(b,w) - {{\rm d} \over {\rm d}t}_{|t=0} {e^{\ad(u(w)+tT(b,w))}-\id \over \ad(u(w)+tT(b,w))}\, T(a,w)\!\!\bigg) \\
  + {\ad(u(w)) e^{\ad(u(w))} \over e^{\ad(u(w))} - \id}\,\big([(e^{-\ad(u(w))}a)_\mfrak{p},(e^{-\ad(u(w))}b)_{\widebar{\mfrak{u}}}]_{\widebar{\mfrak{u}}} + [(e^{-\ad(u(w))}a)_{\widebar{\mfrak{u}}},(e^{-\ad(u(w))}b)_\mfrak{p}]_{\widebar{\mfrak{u}}}\big)
\end{multline*}
multiplied by $\delta(z-w)$. With the help of \eqref{eq:relation II}, we rewrite the expression above into the form
\begin{multline*}
  {\ad(u(w)) \over e^{\ad(u(w))} - \id}\, \bigg[{e^{\ad(u(w))}- \id \over \ad(u(w))}\,T(a,w),{e^{\ad(u(w))}- \id \over \ad(u(w))}\,T(b,w)\bigg] \\
  + {\ad(u(w)) e^{\ad(u(w))} \over e^{\ad(u(w))} - \id}\,\big( [(e^{-\ad(u(w))}a)_\mfrak{p},(e^{-\ad(u(w))}b)_{\widebar{\mfrak{u}}}]_{\widebar{\mfrak{u}}} + [(e^{-\ad(u(w))}a)_{\widebar{\mfrak{u}}},(e^{-\ad(u(w))}b)_\mfrak{p}]_{\widebar{\mfrak{u}}}\big),
\end{multline*}
which finally gives
\begin{multline*}
  {\ad(u(w)) e^{\ad(u(w))}\over e^{\ad(u(w))} - \id}\, \big(e^{-\ad(u(w))}[e^{\ad(u(w))}(e^{-\ad(u(w))}a)_{\widebar{\mfrak{u}}},e^{\ad(u(w))}(e^{-\ad(u(w))}b)_{\widebar{\mfrak{u}}}] \\
  + [(e^{-\ad(u(w))}a)_\mfrak{p},(e^{-\ad(u(w))}b)_{\widebar{\mfrak{u}}}]_{\widebar{\mfrak{u}}} + [(e^{-\ad(u(w))}a)_{\widebar{\mfrak{u}}},(e^{-\ad(u(w))}b)_\mfrak{p}]_{\widebar{\mfrak{u}}} \big) \\
  = {\ad(u(w)) e^{\ad(u(w))}\over e^{\ad(u(w))} - \id}\,[e^{-\ad(u(w))}a, e^{-\ad(u(w))}b]_{\widebar{\mfrak{u}}}  \\
  = {\ad(u(w)) e^{\ad(u(w))}\over e^{\ad(u(w))} - \id}\,(e^{-\ad(u(w))}[a,b])_{\widebar{\mfrak{u}}}.
\end{multline*}
Hence, we get
\begin{align*}
  T([a,b],w)\delta(z-w)= - \big([D(b,z),T(a,w)]- [D(a,z),T(b,w)]\big)
\end{align*}
for all $a,b \in \mfrak{g}$. We are done.}
\vspace{-2mm}


\subsection{Generalized Imaginary Verma modules}

The goal of this section is to construct generalized imaginary Verma modules based on 
the geometric realization introduced and discussed in the previous section.
\medskip

\definition{Let $\mfrak{g}$ be a topological Lie algebra over $\C$ and let $V$ be a topological vector space over $\C$.
A representation $(\sigma,V)$ of $\mfrak{g}$ is called continuous if the corresponding mapping $\mfrak{g} \times V \rarr V$ is continuous.}

Let us consider a continuous $\mfrak{p}_{{\rm nat}}$-module $(\sigma,V)$ such that $\sigma(c)=k\cdot \id_V$ for $k\in \C$. Then the topological vector space $\Pol \Omega_\mcal{K}(\widebar{\mfrak{u}}^*) \otimes_\C\! V$ has a structure of a continuous $\widehat{\mfrak{g}}$-module at level $k$ defined as follows. From Theorem \ref{thm:realization power series}, we have an injective homomorphism
\begin{align}
\pi \colon \widehat{\mfrak{g}} \rarr \eus{R}^{\mfrak{g},\mfrak{p}}_{{\rm loc}}
\end{align}
of topological Lie algebras. Therefore, it is sufficient to show that $\Pol \Omega_\mcal{K}(\widebar{\mfrak{u}}^*) \otimes_\C\! V$ is a continuous $\eus{R}^{\mfrak{g},\mfrak{p}}_{{\rm loc}}$-module.
\medskip

\lemma{\label{lem:a*(z) action}We have
\begin{align}
  \partial^k_za_\alpha^*(z) \bigg(\sum_{n \in \Z} a_n \partial_{x_{\beta,n}}\!\bigg)= -\delta_{\alpha \beta}\,\partial^k_z\bigg(\sum_{n \in \Z} a_n z^n\bigg)
\end{align}
in $\Pol \Omega_\mcal{K}(\widebar{\mfrak{u}}^*)$ for all $k \in \N_0$ and $\alpha,\beta \in \Delta(\mfrak{u})$. Moreover, we have
\begin{align}
  \eus{F}^{\mfrak{g},\mfrak{p}}_{{\rm loc}}(z) \big(\Pol \Omega_\mcal{K}(\widebar{\mfrak{u}}^*)\big) \subset  \Pol \Omega_\mcal{K}(\widebar{\mfrak{u}}^*) \otimes_\C \C(\!(z)\!).  \label{eq:loc function action}
\end{align}
}

\proof{We may write
\begin{align*}
  a_\alpha^*(z) \bigg(\sum_{n \in \Z} a_n \partial_{x_{\beta,n}}\!\bigg) &= \sum_{m \in \Z} \sum_{n  \in \Z} a_n x_{\alpha,m}\partial_{x_{\beta,n}}z^m =
  \sum_{m \in \Z} \sum_{n  \in \Z} a_n [x_{\alpha,m},\partial_{x_{\beta,n}}]z^m \\ & = - \sum_{m \in \Z} \sum_{n  \in \Z} a_n \delta_{\alpha\beta} \delta_{mn}z^m= - \delta_{\alpha\beta}\sum_{n \in \Z} a_n z^n,
\end{align*}
where the second equality follows from the fact that $\Pol \Omega_\mcal{K}(\widebar{\mfrak{u}}^*)$ is an $\eus{A}^{\mfrak{g},\mfrak{p}}$-module.
If we take the derivative of this equation with respect to the formal variable $z$, we obtain the required statement. The rest of the statement is an easy consequence.}

\theorem{Let $(\sigma,V)$ be a continuous $\mfrak{p}_{{\rm nat}}$-module such that $\sigma(c)=k \cdot \id_V$ for $k \in \C$. Then the topological vector space $\Pol \Omega_\mcal{K}(\widebar{\mfrak{u}}^*) \otimes_\C V$ is a continuous $\widehat{\mfrak{g}}$-module at level $k$.}

\proof{Because $V$ is a continuous $\mfrak{p}_{{\rm nat}}$-module and $\Pol \Omega_\mcal{K}(\widebar{\mfrak{u}}^*)$ is a continuous $\eus{F}^{\mfrak{g},\mfrak{p}}$-module, we obtain that the completed tensor product $\Pol \Omega_\mcal{K}(\widebar{\mfrak{u}}^*)\; \smash{\widehat{\otimes}}_\C\, V$ is a continuous $\eus{F}^{\mfrak{g},\mfrak{p}}\; \widehat{\otimes}_\C\; \mfrak{p}_{{\rm nat}}$-module, and therefore also a continuous $\eus{J}^{\mfrak{g},\mfrak{p}}_{{\rm loc}}$-module, since we have $\eus{J}^{\mfrak{g},\mfrak{p}}_{{\rm loc}} \subset \eus{F}^{\mfrak{g},\mfrak{p}}\; \widehat{\otimes}_\C\; \mfrak{p}_{{\rm nat}}$.

In fact, the subspace $\Pol \Omega_\mcal{K}(\widebar{\mfrak{u}}^*) \otimes_\C\! V$ is a $\eus{J}^{\mfrak{g},\mfrak{p}}_{{\rm loc}}$-submodule as follows from the following computation. Let $P(z) \in \eus{F}^{\mfrak{g},\mfrak{p}}_{{\rm loc}}(z)$ and $p \in \Pol \Omega_\mcal{K} (\widebar{\mfrak{u}}^*)$, then from \eqref{eq:loc function action} we have
\begin{align*}
  P(z)p=\sum_{i=1}^r g_i(z)q_i,
\end{align*}
where $q_i \in \Pol \Omega_\mcal{K}(\widebar{\mfrak{u}}^*)$ and $g_i(z) \in \C(\!(z)\!)$. Therefore, we may write
\begin{align*}
  (\Res_{z=0}P(z)h(z)f(z)dz) (p\otimes v) &= \Res_{z=0}(P(z)f(z)p \otimes h(z)v)dz\\
   & = \sum_{i=1}^r \Res_{z=0}(g_i(z)f(z)q_i \otimes h(z)v)dz \\
   & = \sum_{i=1}^r \sum_{n \in \Z} (b_{i,n}q_i \otimes h_nv)
\end{align*}
for $h \in \mfrak{p}$, $f(z) \in \C(\!(z)\!)$, $p \in \Pol \Omega_\mcal{K}(\widebar{\mfrak{u}}^*)$ and $v \in V$, where $g_i(z)f(z)=\sum_{n\in \Z} b_{i,n}z^n$ for $i=1,2,\dots,r$. Then for all $i=1,2,\dots,r$, we have
\begin{align*}
\sum_{n \in \Z} (b_{i,n}q_i \otimes h_nv) &= \sum_{n=n_i}^\infty (b_{i,n}q_i \otimes h_nv) = \lim_{N \rarr \infty} \sum_{n=n_i}^N  (b_{i,n}q_i \otimes h_nv) \\
&= \lim_{N \rarr \infty} \sum_{n=n_i}^N  (q_i \otimes b_{i,n}h_nv) = \lim_{N \rarr \infty}\! \bigg(q_i \otimes \sum_{n=n_i}^N b_{i,n}h_nv\!\bigg) \\
&= q_i \otimes \sum_{n \in \Z}b_{i,n}h_nv,
\end{align*}
since $b_{i,n} = 0$ for $n < n_i$ and $V$ is a continuous $\mfrak{p}_{{\rm nat}}$-module. Hence, the subspace $\Pol \Omega_\mcal{K}(\widebar{\mfrak{u}}^*) \otimes_\C\! V$ is a $\eus{J}^{\mfrak{g},\mfrak{p}}_{{\rm loc}}$-submodule. Because $\Pol \Omega_\mcal{K}(\widebar{\mfrak{u}}^*)$ is a continuous $\eus{A}^{\mfrak{g},\mfrak{p},{\rm op}}$-module, we get also a structure of a continuous $\eus{A}^{\mfrak{g},\mfrak{p},{\rm op}}$-module on $\Pol \Omega_\mcal{K}(\widebar{\mfrak{u}}^*) \otimes_\C\! V$. Therefore $\Pol \Omega_\mcal{K}(\widebar{\mfrak{u}}^*) \otimes_\C\! V$ is a continuous $\smash{\eus{A}^{\mfrak{g},\mfrak{p},{\rm op}}_{\leq 1, {\rm loc}}}$-module.

However, since the Lie algebra $\eus{R}^{\mfrak{g},\mfrak{p}}_{{\rm loc}}$ is not a direct sum (only semidirect sum) of the Lie algebras $\smash{\eus{A}^{\mfrak{g},\mfrak{p},{\rm op}}_{\leq 1, {\rm loc}}}$ and $\eus{J}^{\mfrak{g},\mfrak{p}}_{{\rm loc}}$, we need to verify that $\Pol \Omega_\mcal{K}(\widebar{\mfrak{u}}^*) \otimes_\C\! V$ is an $\eus{R}^{\mfrak{g},\mfrak{p}}_{{\rm loc}}$-module. For $P(z) \in \eus{F}^{\mfrak{g},\mfrak{p}}_{{\rm loc}}(z)$, $h \in \mfrak{p}$ and $Q(w)=T(w)S(w)$, where $T(w)$ is a differential polynomial in $a_\alpha(w)$ for $\alpha \in \Delta(\mfrak{u})$ of degree at most $1$ and $S(w) \in \eus{F}^{\mfrak{g},\mfrak{p}}_{{\rm loc}}(w)$, we have
\begin{multline*}
  (\Res_{w=0}Q(w)g(w)dw\Res_{z=0}P(z)h(z)f(z)dz)(p\otimes v)\\
  =\Res_{z=0,w=0}(f(z)g(w)Q(w)P(z)p  \otimes h(z)v)dzdw
\end{multline*}
and
\begin{multline*}
  (\Res_{z=0}P(z)h(z)f(z)dz\Res_{w=0}Q(w)g(w)dw)(p\otimes v)\\
  =\Res_{z=0,w=0}(f(z)g(w)P(z)Q(w)p  \otimes h(z)v)dzdw
\end{multline*}
for all $p \in \Pol \Omega_\mcal{K} (\widebar{\mfrak{u}}^*)$ and $v \in V$.
Taking the difference of the previous two equations, we obtain for the right hand side
\begin{multline*}
  \Res_{z=0,w=0}([Q(w)g(w),P(z)f(z)]p \otimes h(z)v)dzdw\\
  =(\Res_{z=0,w=0}[Q(w)g(w),P(z)f(z)]h(z)dzdw)(p \otimes v) \\
  =[\Res_{w=0}Q(w)g(w)dw,\Res_{z=0}P(z)h(z)f(z)dz](p \otimes v),
\end{multline*}
which we required. A similar computation can be done, if we replace $h(z)$ by the central element $c$. Hence, the commutation relations between $\smash{\eus{A}^{\mfrak{g},\mfrak{p},{\rm op}}_{\leq 1, {\rm loc}}}$ and $\eus{J}^{\mfrak{g},\mfrak{p}}_{{\rm loc}}$ are satisfied in $\Pol \Omega_\mcal{K}(\widebar{\mfrak{u}}^*) \otimes_\C\! V$. The continuity follows immediately.}

Our second main result is the following theorem, which relates generalized imaginary Verma modules $\mathbb{M}_{\sigma,k,\mfrak{p}}(V)$ and $\Pol \Omega_\mcal{K} (\widebar{\mfrak{u}}^*) \otimes_\C\! V$.
\medskip

\theorem{\label{thm:isomorphism}
Let $(\sigma,V)$ be a continuous $\mfrak{p}_{{\rm nat}}$-module such that $\sigma(c)=k\cdot \id_V$ for $k \in \C$. Then we have
\begin{align}
  \mathbb{M}_{\sigma,k,\mfrak{p}}(V) \simeq \Pol \Omega_\mcal{K} (\widebar{\mfrak{u}}^*) \otimes_\C\! V
\end{align}
as continuous $\widehat{\mfrak{g}}$-modules.}

\proof{From the previous theorem we obtained a continuous $\widehat{\mfrak{g}}$-module $\Pol \Omega_\mcal{K} (\widebar{\mfrak{u}}^*) \otimes_\C\! V$ at level $k$ for any continuous $\mfrak{p}_{{\rm nat}}$-module $(\sigma,V)$ such that $\sigma(c)=k \cdot \id_V$ for $k\in \C$. We shall show that this $\widehat{\mfrak{g}}$-module is the generalized imaginary Verma module $\mathbb{M}_{\sigma,k,\mfrak{p}}(V)$. From Theorem \ref{thm:realization power series} we have that the canonical mapping
\begin{align*}
  \varphi_0 \colon V \rarr \Pol \Omega_\mcal{K} (\widebar{\mfrak{u}}^*) \otimes_\C\! V
\end{align*}
is a homomorphism of continuous $\mfrak{p}_{{\rm nat}}$-modules and gives rise to a homomorphism
\begin{align*}
  \varphi \colon \mathbb{M}_{\sigma,k,\mfrak{p}}(V) \rarr \Pol \Omega_\mcal{K} (\widebar{\mfrak{u}}^*) \otimes_\C\! V
\end{align*}
of continuous $\widehat{\mfrak{g}}$-modules. Since $U(\widebar{\mfrak{u}}_{{\rm nat}})$ and $S(\widebar{\mfrak{u}}_{{\rm nat}})\simeq \Pol \Omega_\mcal{K} (\widebar{\mfrak{u}}^*)$ have natural increasing filtrations, we get increasing filtrations on $U(\widebar{\mfrak{u}}_{{\rm nat}}) \otimes_\C\! V$ and $S(\widebar{\mfrak{u}}_{{\rm nat}}) \otimes_\C\! V$. Moreover, we have that $\mathbb{M}_{\sigma,k,\mfrak{p}}(V)$ and $\Pol \Omega_\mcal{K} (\widebar{\mfrak{u}}^*) \otimes_\C\! V$ are filtered $\widebar{\mfrak{u}}_{{\rm nat}}$-modules and $\varphi$ is a homomorphism of filtered $\widebar{\mfrak{u}}_{{\rm nat}}$-modules. To prove that $\varphi$ is an isomorphism, it is enough to show that the associated mapping
\begin{align*}
  \gr \varphi \colon \gr \mathbb{M}_{\sigma,k,\mfrak{p}}(V) \rarr \gr \Pol \Omega_\mcal{K} (\widebar{\mfrak{u}}^*) \otimes_\C\!V
\end{align*}
of graded $\widebar{\mfrak{u}}_{{\rm nat}}$-modules is an isomorphism. Let us define
\begin{align*}
  f_\alpha^g = \Res_{z=0}f_\alpha(z)g(z)dz\qquad \text{and} \qquad a_\alpha^g= \Res_{z=0}a_\alpha(z)g(z)dz
\end{align*}
for $g \in \C(\!(z)\!)$ and $\alpha \in \Delta(\mfrak{u})$. Then the set $\{f_\alpha^g;\, g \in \C(\!(z)\!),\, \alpha\in \Delta(\mfrak{u})\}$ generates $U(\widebar{\mfrak{u}}_{{\rm nat}})$ and the set $\{a_\alpha^g;\, g \in \C(\!(z)\!),\, \alpha\in \Delta(\mfrak{u})\}$ generates $S(\widebar{\mfrak{u}}_{{\rm nat}})$.
Further, from \eqref{eq:action power series op nilradical} we have
\begin{align*}
  \pi(f(z)) =-\sum_{\beta\in \Delta(\mfrak{u})} \sum_{k=0}^\infty B_k a_\beta(z)[(\ad(u(z)))^k(f)]_\beta,
\end{align*}
where the Bernoulli numbers $B_k$ are determined by the generating series
\begin{align*}
  {x \over e^x-1}= \sum_{k=0}^\infty B_k {x^k \over k!}
\end{align*}
for $0 \neq x \in \R$. Furthermore, if we denote
\begin{align*}
  R(f_\alpha,z)= -\sum_{\beta\in \Delta(\mfrak{u})} \sum_{k=1}^\infty B_k a_\beta(z)[(\ad(u(z)))^k(f_\alpha)]_\beta,
\end{align*}
then we may write
\begin{align*}
  \pi(f_\alpha(z))= -a_\alpha(z) + R(f_\alpha,z).
\end{align*}
Moreover, we have
\begin{align*}
  (\Res_{z=0} R(f_\alpha,z)g(z)dz)(S(\widebar{\mfrak{u}}_{{\rm nat}})_m) \subset S(\widebar{\mfrak{u}}_{{\rm nat}})_m
\end{align*}
for all $m\in \N_0$. Since
\begin{align*}
  \pi(f_\alpha^g)=-a_\alpha^g + \Res_{z=0} R(f_\alpha,z)g(z)dz,
\end{align*}
we may write
\begin{align*}
  \begin{aligned}
    (\gr \varphi)_m(f^{g_1}_{\alpha_1}\dots f^{g_m}_{\alpha_m}\otimes v) &= \varphi(f^{g_1}_{\alpha_1}\dots f^{g_m}_{\alpha_m}\otimes v)\ {\rm mod}\ S(\widebar{\mfrak{u}}_{{\rm nat}})_{m-1} \otimes_\C\! V \\
    &=\pi(f^{g_1}_{\alpha_1})\dots \pi(f^{g_m}_{\alpha_m})\varphi_0(v)\ {\rm mod}\ S(\widebar{\mfrak{u}}_{{\rm nat}})_{m-1} \otimes_\C\! V \\
    &=(-1)^m a_{\alpha_1}^{g_1} \dots a_{\alpha_m}^{g_m}\otimes v \ {\rm mod}\ S(\widebar{\mfrak{u}}_{{\rm nat}})_{m-1} \otimes_\C\! V
  \end{aligned}
\end{align*}
for all $m\in \N_0$, $g_1,g_2,\dots,g_m \in \C(\!(z)\!)$ and $\alpha_1,\alpha_2,\dots,\alpha_m \in \Delta(\mfrak{u})$, where we used $\varphi_0(v) \in S(\widebar{\mfrak{u}}_{{\rm nat}})_0 \otimes_\C\! V$ in the last equality. This gives that $\gr \varphi$ is an isomorphism and so $\varphi$ is also an isomorphism.}

In order to employ the irreducibility criterion given in 
\cite{Futorny-Kashuba2016}, we restrict to the class $(\sigma,V)$ 
of weight tensor inducing modules with non-zero central 
charge. The definition of the notion of tensor module  
can be found in \cite{Futorny-Kashuba2016}, conjecturally 
all modules belong to this class.
\medskip

\corollary{Let $(\sigma,V)$ be a weight tensor continuous irreducible $\mfrak{l}_{{\rm nat}}$-module regarded as $\mfrak{p}_{{\rm nat}}$-module with the trivial action of $\mfrak{u}_{{\rm nat}}$ such that $\sigma(c)=k\cdot \id_V$ for $k \in \C$.
If $k\neq 0$, then $\Pol \Omega_\mcal{K} (\widebar{\mfrak{u}}^*) \otimes_\C\! V$ is a continuous irreducible $\widehat{\mfrak{g}}$-module.}

\proof{By \cite{Futorny-Kashuba2016} and the assumptions of the corollary is the generalized imaginary Verma module $\mathbb{M}_{\sigma,k,\mfrak{p}}(V)$ irreducible.}
\vspace{-2mm}


\section{Generalized Imaginary Verma modules for $\smash{\widehat{\mfrak{sl}}(n+1,\C)}$}
\label{sec4}

In the present section we apply the previous general exposition to the case of generalized imaginary Verma modules for the pair given by the complex simple Lie algebra $\mfrak{sl}(n+1,\C)$, $n\in \N$, and its maximal parabolic subalgebra with the commutative nilradical. In the Dynkin diagrammatic notation, this type of parabolic subalgebra is determined by omitting the first simple root in the Dynkin diagram for $\mfrak{sl}(n+1,\C)$.


\subsection{Representation theoretical conventions}

In the rest of the section we consider the complex simple Lie algebra $\mfrak{g}=\mfrak{sl}(n+1,\C)$, $n \in \N$. The Cartan subalgebra $\mfrak{h}$ of $\mfrak{g}$ is given by diagonal
matrices
\begin{align}
  \mfrak{h}=\{\diag(a_1,a_2,\dots,a_{n+1});\, a_1,a_2,\dots,a_{n+1} \in \C,\ {\textstyle \sum_{i=1}^{n+1}}a_i=0\}.
\end{align}
For $i=1,2,\dots,n+1$ we define $\veps_i \in \mfrak{h}^*$ by $\veps_i(\diag(a_1,a_2,\dots,a_{n+1}))=a_i$. Then the root system of $\mfrak{g}$ with respect to $\mfrak{h}$ is
$\Delta = \{\veps_i-\veps_j;\, 1 \leq i \neq j \leq n+1\}$. The root space
 $\mfrak{g}_{\veps_i-\veps_j}$ is the complex linear span of $E_{ij}$,
the $(n+1 \times n+1)$-matrix having $1$ at the intersection of the $i$-th row and $j$-th column and $0$ elsewhere. The positive root system is $\Delta^+=\{\veps_i-\veps_j;\, 1\leq i<j \leq n+1\}$ with the set of simple roots $\Pi=\{\alpha_1,\alpha_2,\dots,\alpha_n\}$, $\alpha_i=\veps_i-\veps_{i+1}$ for $i=1,2,\dots,n$, and the fundamental weights are
$\omega_i= \smash{\sum_{j=1}^i \veps_j}$ for $i=1,2,\dots,n$.
The solvable Lie subalgebras $\mfrak{b}$ and $\widebar{\mfrak{b}}$ defined as the direct sum of positive and negative root
spaces together with the Cartan subalgebra are called the standard Borel subalgebra and the opposite standard Borel subalgebra of $\mfrak{g}$, respectively.
The subset $\Sigma=\{\alpha_2,\alpha_3,\dots,\alpha_n\}$
of $\Pi$ generates the root subsystem $\Delta_\Sigma$ in $\mfrak{h}^*$, and we associate to $\Sigma$ the standard parabolic
subalgebra $\mfrak{p}$ of $\mfrak{g}$ by $\mfrak{p} = \mfrak{l} \oplus \mfrak{u}$.
The reductive Levi subalgebra $\mfrak{l}$ of $\mfrak{p}$ is defined through
\begin{align}
  \mfrak{l}= \mfrak{h} \oplus \bigoplus_{\alpha \in \Delta_\Sigma} \mfrak{g}_\alpha,
\end{align}
and the nilradical $\mfrak{u}$ of $\mfrak{p}$ and the opposite nilradical $\widebar{\mfrak{u}}$ are given by
\begin{align}
  \mfrak{u}= \bigoplus_{\alpha \in \Delta^+ \setminus \Delta_\Sigma^+}\mfrak{g}_\alpha \qquad \text{and}
  \qquad \widebar{\mfrak{u}}= \bigoplus_{\alpha \in \Delta^+ \setminus \Delta_\Sigma^+} \mfrak{g}_{-\alpha},
\end{align}
respectively. We define the $\Sigma$-height $\htt_\Sigma(\alpha)$ of $\alpha \in \Delta$ by
\begin{align}
  \htt_\Sigma\!\big({\textstyle \sum_{i=1}^n} a_i \alpha_i\big) = a_1,
\end{align}
so $\mfrak{g}$ is a $|1|$-graded Lie algebra with respect to the grading given by
$\mfrak{g}_i = \bigoplus_{\alpha \in \Delta,\, \htt_\Sigma(\alpha)=i} \mfrak{g}_\alpha$
for $0 \neq i \in \Z$, and
$\mfrak{g}_0= \mfrak{h} \oplus \smash{\bigoplus}_{\alpha \in \Delta,\, \htt_\Sigma(\alpha)=0} \mfrak{g}_\alpha$. Moreover, we have
$\mfrak{u}=\mfrak{g}_1$, $\widebar{\mfrak{u}}=\mfrak{g}_{-1}$ and $\mfrak{l}= \mfrak{g}_0$.

Furthermore, the normalized $\mfrak{g}$-invariant symmetric bilinear form $(\cdot\,,\cdot)\colon \mfrak{g} \otimes_\C \mfrak{g} \rarr \C$ on $\mfrak{g}$ is given by
\begin{align}
  (a,b) = \tr(ab)
\end{align}
for $a,b \in \mfrak{g}$, since we have $\theta=\veps_1-\veps_{n+1}$, and the dual Coxeter number is $h^\vee=n+1$.

Now, let $\{f_1,f_2,\dots,f_n\}$ be a basis of the root spaces in the opposite nilradical $\widebar{\mfrak{u}}$ given by
\begin{align}
  f_i =\begin{pmatrix}
    0 & 0  \\
    1_i & 0
  \end{pmatrix}
\end{align}
and let $\{e_1,e_2\dots,e_n\}$ be a basis of the root spaces in the nilradical $\mfrak{u}$ defined by
\begin{align}
  e_i =\begin{pmatrix}
    0 & 1_i^{\rm T} \\
    0 & 0
  \end{pmatrix}\!.
\end{align}
The Levi subalgebra $\mfrak{l}$ of $\mfrak{p}$ is the linear span of
\begin{align}
  h =\begin{pmatrix}
    1 & 0  \\
    0 & -{1\over n} I_n  \\
  \end{pmatrix} \qquad \text{and} \qquad
  h_A =\begin{pmatrix}
    0 & 0 \\
    0 & A
  \end{pmatrix}\!,
\end{align}
where $A \in M_{n\times n}(\C)$ satisfies $\tr A = 0$. Moreover, the element $h$ forms a basis of the center $\mfrak{z}(\mfrak{l})$ of the Levi subalgebra $\mfrak{l}$.


\subsection{The embedding of $\widehat{\mfrak{g}}$ into the Lie algebra $\eus{R}^{\mfrak{g},\mfrak{p}}_{{\rm loc}}$}

Let us denote by $\{x_i;\, i=1,2,\dots,n\}$ the linear coordinate functions on $\widebar{\mfrak{u}}$ with respect to the basis $\{f_i;\, i=1,2,\dots,n\}$ of the opposite nilradical $\widebar{\mfrak{u}}$. Then the set $\{f_i \otimes t^m;\, m\in \Z,\, i=1,2,\dots,n\}$ forms a topological basis of $\mcal{K}(\widebar{\mfrak{u}})=\widebar{\mfrak{u}}_{{\rm nat}}$, and the set $\{x_i \otimes t^{-m-1}dt;\, m\in\Z,\, i=1,2,\dots,n\}$ forms a dual topological basis of $\Omega_\mcal{K}(\widebar{\mfrak{u}}^*)\simeq (\widebar{\mfrak{u}}_{{\rm nat}})^*$ with respect to the pairing \eqref{eq:non-degenerate form}. We denote $x_{i,m} = x_i \otimes t^{-m-1}dt$ and $\partial_{x_{i,m}}=f_i \otimes t^m$ for $m\in \Z$ and $i=1,2,\dots,n$. The Weyl algebra $\smash{\eus{A}_{\mcal{K}(\widebar{\mfrak{u}})}}$ is topologically generated by $\{x_{i,m},\partial_{x_{i,m}};\, m \in \Z,\, i=1,2,\dots,n\}$ with the canonical commutation relations. Furthermore, we define the formal distributions $a_i(z), a^*_i(z) \in \smash{\eus{A}_{\mcal{K}(\widebar{\mfrak{u}})}[[z^{\pm 1}]]}$ by
\begin{align}
  a_i(z)= \sum_{m\in \Z} a_{i,m} z^{-m-1}  \qquad \text{and} \qquad a^*_i(z)= \sum_{m\in \Z}  a^*_{i,m}\, z^{-m},
\end{align}
where $a_{i,m}=\partial_{x_{i,m}}$ and $a^*_{i,m}=x_{i,-m}$, for $i=1,2,\dots,n$. Finally, let us introduce the formal power series $u(z) \in \widebar{\mfrak{u}} \otimes_\C \smash{\eus{P}^{\mfrak{g},\mfrak{p}}_{{\rm loc}}(z)}$ by
\begin{align}
  u(z)=\sum_{i=1}^n a_i^*(z)f_i.
\end{align}
\smallskip

\theorem{\label{thm:imaginary sl realization series}
The embedding of $\widehat{\mfrak{g}}$ into $\eus{R}^{\mfrak{g},\mfrak{p}}_{{\rm loc}}$ is given by
\begin{enumerate}
  \item[1)]
  \begin{align}
    \pi(f_i(z))=-a_i(z) \label{eq:realization sl nilradical opposite}
  \end{align}
  for $i=1,2,\dots,n$;
  \item[2)]
  \begin{align} \label{eq:realization sl levi}
    \begin{aligned}
      \pi(c)&=c,\\
      \pi(h(z))&=\big(1+{\textstyle {1\over n}}\big){\textstyle \sum_{j=1}^n} a_j(z)a^*_j(z)+h(z), \\
      \pi(h_A(z))&=-{\textstyle \sum_{r,s=1}^n} a_{r,s} a_r(z)a^*_s(z)+h_A(z)
    \end{aligned}
  \end{align}
  for all $A=(a_{r,s}) \in M_{n \times n}(\C)$ satisfying $\tr A=0$;
  \item[3)]
  \begin{align} \label{eq:realization sl nilradical}
  \begin{aligned}
    \pi(e_i(z))&={\textstyle \sum_{j=1}^n} a_j(z)a^*_j(z)a^*_i(z)-\partial_z a^*_i(z)c \\ &\quad +e_i(z) +a^*_i(z)h(z)- {\textstyle \sum_{j=1}^n} a^*_j(z)h_{E_{ji}-{1 \over n}I_n \delta_{ij}}(z)
  \end{aligned}
  \end{align}
  for $i=1,2,\dots,n$.
\end{enumerate}}

\proof{By Theorem \ref{thm:realization power series} we have
\begin{align*}
  \pi(a(z))=-\sum_{i=1}^n a_i(z)\bigg[{\ad(u(z)) \over e^{\ad(u(z))}-\id}\,a\bigg]_i
\end{align*}
for $a \in \widebar{\mfrak{u}}$. Because $\ad(u(z))(a)=0$, the expansion of the formal power series in $\ad(u(z))$ implies
\begin{align*}
  \pi(a(z))=-\sum_{i=1}^n a_i(z)[a]_i,
\end{align*}
which gives \eqref{eq:realization sl nilradical opposite}. Similarly, from Theorem \ref{thm:realization power series} we get
\begin{align*}
  \pi(a(z))=  \sum_{i=1}^n a_i(z) [\ad(u(z))(a)]_i + a(z)
\end{align*}
for $a \in \mfrak{l}$, and we obtain \eqref{eq:realization sl levi}. Finally, from Theorem \ref{thm:realization power series} we have
\begin{multline*}
  \pi(a(z))= - \sum_{i=1}^n a_i(z) \bigg[{\ad(u(z)) e^{\ad(u(z))} \over e^{\ad(u(z))}- \id}\, (e^{-\ad(u(z))}a)_{\widebar{\mfrak{u}}}\bigg]_i \\ + (e^{-\ad(u(z))} a(z))_\mfrak{p} - \bigg({e^{\ad(u(z))}-\id \over \ad(u(z))}\, \partial_z u(z),a\!\bigg)c,
\end{multline*}
for $a \in \mfrak{u}$. Since $\mfrak{g}$ is $|1|$-graded, we have $(\ad(u(z)))^3(a)=0$. Therefore, we get
\begin{align*}
  (e^{-\ad(u(z))}a)_{\widebar{\mfrak{u}}}={\textstyle {1 \over 2}}(\ad(u(z)))^2(a) \qquad \text{and} \qquad (a^{-\ad(u(z))}a(z))_\mfrak{p}=a(z)-\ad(u(z))(a(z)).
\end{align*}
Hence, we may write
\begin{align*}
  \pi(a(z))= - {1\over 2}\sum_{i=1}^n a_i(z)[(\ad(u(z)))^2(a)]_i  + a(z)-\ad(u(z))(a(z)) - (\partial_zu(z),a)c,
\end{align*}
and \eqref{eq:realization sl nilradical} follows from the commutation relations in $\mfrak{g}$. This completes the proof.}

We can write explicitly the action of the topological generators of $\widehat{\mfrak{g}}$.
\medskip

\theorem{\label{thm:imaginary sl realization}
The embedding of $\widehat{\mfrak{g}}$ into $\eus{R}^{\mfrak{g},\mfrak{p}}_{{\rm loc}}$ is given by
\begin{enumerate}
  \item[1)]
  \begin{align}
      \pi(f_{i,m})=-\partial_{x_{i,m}}
  \end{align}
  for $i=1,2,\dots,n$ and $m\in \Z$;
  \item[2)]
  \begin{align}
    \begin{aligned}
      \pi(c) &= c, \\
      \pi(h_m) & = \big(1+{\textstyle {1 \over n}}\big){\textstyle \sum_{j=1}^n \sum_{k \in \Z}}\, \partial_{x_{j,k+m}} x_{j,k}+ h_m, \\
      \pi(h_{A,m}) &= -{\textstyle \sum_{r,s=1}^n \sum_{k\in \Z}}\, a_{r,s}\partial_{x_{r,k+m}} x_{s,k} + h_{A,m}
    \end{aligned}
  \end{align}
  for all $A=(a_{r,s}) \in M_{n \times n}(\C)$ satisfying $\tr A=0$ and $m\in \Z$;
  \item[3)]
  \begin{align}
  \begin{aligned}
    \pi(e_{i,m})&= {\textstyle \sum_{j=1}^n \sum_{k,\ell \in \Z}} \,\partial_{x_{j,k+\ell+m}}x_{i,k}x_{j,\ell}+mx_{i,-m}c \\ &\quad  + e_{i,m}+{\textstyle \sum_{k\in \Z}}\, x_{i,k} h_{k+m} - {\textstyle \sum_{j=1}^n \sum_{k\in \Z}}\, x_{j,k} h_{E_{ji}-{1 \over n}I_n \delta_{ij},k+m}
  \end{aligned}
  \end{align}
  for $i=1,2,\dots,n$ and $m\in \Z$.
\end{enumerate}}

\proof{The proof follows easily from Theorem \ref{thm:imaginary sl realization series}, if we expand the corresponding formal power series.}

Let $\sigma \colon \mfrak{p}_{{\rm nat}} \rarr \mfrak{gl}(V)$ be a continuous $\mfrak{p}_{{\rm nat}}$-module such that $\sigma(c)=k \cdot \id_V$ for $k\in \C$. The generalized imaginary Verma module at level $k$ is
\begin{align}
\mathbb{M}_{\sigma,k,\mfrak{p}}(V) = U(\widehat{\mfrak{g}}) \otimes_{U(\mfrak{p}_{{\rm nat}})}\! V \simeq \Pol \Omega_\mcal{K}(\widebar{\mfrak{u}}^*) \otimes_\C V,
\end{align}
where the corresponding action of $\widehat{\mfrak{g}}$ on $\Pol \Omega_\mcal{K}(\widebar{\mfrak{u}}^*) \otimes_\C V$ is given through the mapping
\begin{align}
 \pi \colon \widehat{\mfrak{g}} \rarr \eus{R}^{\mfrak{g},\mfrak{p}}_{{\rm loc}}
\end{align}
in Theorem \ref{thm:imaginary sl realization}. 


\subsection{Imaginary Verma modules for $\smash{\widehat{\mfrak{sl}}(2,\C)}$}

For the reader's convenience, in this subsection we describe explicitly the imaginary Verma modules for the affine Kac-Moody algebra $\smash{\widehat{\mfrak{sl}}(2,\C)}$. Despite the fact that it corresponds to the choice $n=1$ in the previous subsection, it is convenient to state it separately, as the simplest possible case.

Let us consider the standard basis of the Lie algebra $\mfrak{g}=\mfrak{sl}(2,\C)$ given by
\begin{align}
  e= \begin{pmatrix}
    0 & 1 \\
    0 & 0
  \end{pmatrix}\!,
  \qquad
  h = \begin{pmatrix}
    1 & 0 \\
    0 & -1
  \end{pmatrix}\!,
  \qquad
  f = \begin{pmatrix}
    0 & 0 \\
    1 & 0
  \end{pmatrix}\!.
\end{align}
The normalized $\mfrak{g}$-invariant symmetric bilinear form $(\cdot\,,\cdot)\colon \mfrak{g} \otimes_\C \mfrak{g} \rarr \C$ on $\mfrak{g}$ is given by $(a,b) = \tr(ab)$ for $a,b \in \mfrak{g}$, and the dual Coxeter number is $h^\vee=2$. In particular, we have
\begin{align}
  (e,f)=1, \qquad (h,h)=2, \qquad (f,e)=1,
\end{align}
and zero otherwise.

We consider the Borel subalgebra $\mfrak{b}=\C e\oplus \C h$ with the corresponding nilradical $\mfrak{n}=\C e$ and opposite nilradical $\widebar{\mfrak{n}}=\C f$. Let $x \colon \widebar{\mfrak{n}} \rarr \C$ be a linear coordinate function on $\widebar{\mfrak{n}}$ defined by $x(f)=1$. Then the set $\{f \otimes t^n;\, n\in\Z\}$ forms a topological basis of $\mcal{K}(\widebar{\mfrak{n}})=\widebar{\mfrak{n}}_{{\rm nat}}$, and the set $\{x \otimes t^{-n-1}dt;\, n \in \Z\}$ forms a dual topological basis of $\Omega_\mcal{K}(\widebar{\mfrak{n}}^*) \simeq (\widebar{\mfrak{n}}_{{\rm nat}})^*$ with respect to the pairing \eqref{eq:non-degenerate form}. We define the formal distributions $a(z), a^*(z) \in \eus{A}_{\mcal{K}(\widebar{\mfrak{n}})}[[z^{\pm 1}]]$ by
\begin{align}
  a(z)=  \sum_{n\in \Z} a_n z^{-n-1} = \sum_{n \in \Z} \partial_{x_n} z^{-n-1} \qquad \text{and} \qquad a^*(z)= \sum_{n\in \Z} a^*_n z^{-n} = \sum_{n \in \Z} x_{-n}z^{-n},
\end{align}
so that we have $a(z)=a_\alpha(z)$ and $a^*(z)=a^*_\alpha(z)$. Finally, we introduce the formal power series $u(z) \in \widebar{\mfrak{n}} \otimes_\C \smash{\eus{P}^{\mfrak{g},\mfrak{b}}_{{\rm loc}}(z)}$ by
\begin{align}
  u(z)=a^*(z)f.
\end{align}
\smallskip

\theorem{\label{thm:imaginary sl2 Verma}
The embedding of $\widehat{\mfrak{g}}$ into $\eus{R}^{\mfrak{g},\mfrak{b}}_{{\rm loc}}$ is given by
\begin{align}
  \begin{aligned}
      \pi(c) &= c, \\
      \pi(f_n)&=-\partial_{x_n}, \\
      \pi(h_n) & = 2\,{\textstyle \sum_{k \in \Z}}\, \partial_{x_{k+n}}x_k+ h_n, \\
      \pi(e_n)&=\smash{{\textstyle \sum_{k,\ell \in \Z}}\, \partial_{x_{k+\ell+n}}x_kx_\ell+nx_{-n}c +{\textstyle \sum_{k\in \Z}}\, x_k h_{k+n} + e_n}
  \end{aligned}
\end{align}
for $n\in \Z$, or equivalently by
\begin{align}
  \begin{aligned}
    \pi(c)&=c, \\
    \pi(f(z))&=-a(z), \\
    \pi(h(z))&=2a(z)a^*(z)+h(z), \\
    \pi(e(z))&=a(z)a^*(z)^2-\partial_z a^*(z)c+a^*(z)h(z)+e(z),
  \end{aligned}
  \end{align}
if we use the formal distributions.}

\proof{The proof is a consequence of Theorem \ref{thm:imaginary sl realization} and Theorem \ref{thm:imaginary sl realization series} for $n=1$.}

Let $\sigma \colon \mfrak{b}_{{\rm nat}} \rarr \mfrak{gl}(V)$ be a continuous $\mfrak{b}_{{\rm nat}}$-module such that $\sigma(c)=k\cdot \id_V$ for $k\in \C$. The imaginary Verma module at level $k$ is
\begin{align}
\mathbb{M}_{\sigma,k,\mfrak{b}}(V) = U(\widehat{\mfrak{g}}) \otimes_{U(\mfrak{p}_{{\rm nat}})}\! V \simeq \Pol \Omega_\mcal{K}(\widebar{\mfrak{n}}^*) \otimes_\C V,
\end{align}
where the corresponding action of $\widehat{\mfrak{g}}$ on $\Pol \Omega_\mcal{K}(\widebar{\mfrak{n}}^*) \otimes_\C V$ is given through the mapping
\begin{align}
\pi \colon \widehat{\mfrak{g}} \rarr \eus{R}^{\mfrak{g},\mfrak{b}}_{{\rm loc}}
\end{align}
in Theorem \ref{thm:imaginary sl2 Verma}. 


\section*{Acknowledgments}
V.\,Futorny is supported in part by the CNPq grant (301320/2013-6) and by the
Fapesp grant (2014/09310-5). He gratefully acknowledges the hospitality and excellent
working conditions at the  Charles University where this work was done. L.\,Křižka is supported by PRVOUK p47. He is grateful to the University of Sa\~{o} Paulo for the hospitality, where part of this work was done. P.\,Somberg acknowledges the financial support from the grant GA\,P201/12/G028.



\providecommand{\bysame}{\leavevmode\hbox to3em{\hrulefill}\thinspace}
\providecommand{\MR}{\relax\ifhmode\unskip\space\fi MR }
\providecommand{\MRhref}[2]{%
  \href{http://www.ams.org/mathscinet-getitem?mr=#1}{#2}
}
\providecommand{\href}[2]{#2}

\end{document}